\documentclass[notitlepage,leqno,10pt]{amsart}
\usepackage{amsfonts,amssymb,cite}
\usepackage[latin1]{inputenc}
\usepackage{enumerate}
\usepackage{amssymb}
\catcode`\@=11 \@addtoreset{equation}{section}

\catcode`\@=12
\usepackage{latexsym}
\usepackage{textcomp}
\usepackage{amsmath}
\usepackage{amsfonts}
\usepackage{amssymb}
\usepackage{mathrsfs}

\renewcommand{\a }{\alpha}
\renewcommand{\b }{\beta}
\renewcommand{\d}{\delta}

\newcommand{\D }{\Delta }

\newcommand{\e }{\varepsilon }
\newcommand{\g }{\gamma}

\newcommand{\G }{\Gamma}
\renewcommand{\l }{\lambda }

\newcommand{\n }{\nabla }

\newcommand{\s }{\sigma }
\newcommand{\Sig }{\Sigma}
\renewcommand{\t }{\tau }

\renewcommand{\o }{\omega }
\renewcommand{\O }{\Omega }

\newcommand{\ov}{\overline}

\newcommand{\be}{\begin{equation}}
\newcommand{\ee}{\end{equation}}

\newcommand{\R}{\mathbb{R}}
\newcommand{\N}{\mathbb{N}}

\newcommand{\de}{\partial}

\newcommand{\ti}{\tilde}

\newcommand{\M}{\mathcal{M}}

\newcommand{\ra}{{\rangle}}
\newcommand{\la}{{\langle}}

\newcommand{\calV }{\mathcal{V}}

\newcommand{\calG }{\mathcal{G}}

\newcommand{\calL }{\mathcal{L}}

\newcommand{\calF }{\mathcal{F}}

\newcommand{\calM}{{\mathcal M}}
\newcommand{\calT}{{\mathcal T}}

\newcommand{\mbL }{\mathbb{L}}

\newtheorem{Theorem}{Theorem}[section]
\newtheorem{Lemma}[Theorem]{Lemma}
\newtheorem{Proposition}[Theorem]{Proposition}

\newtheorem{Remark}[Theorem]{Remark}
\newtheorem{Definition}[Theorem]{Definition}

\def\proof{\noindent{{\bf Proof. }}}
\def\square{\vbox{
\hrule height .4pt \hbox{\vrule width .4pt height 7pt \kern 7pt
\vrule width .4pt} \hrule height .4pt }}
\def\QED{\hfill {$\square$}\goodbreak \medskip}
\def\R{{\mathbb R}}

\font\sc=cmcsc9  \textwidth=14truecm
\begin{document}
\title[Over-determined Problems ]{ Serrin's over-determined Problem on Riemannian
Manifolds}
\author{Mouhamed Moustapha Fall}
\author{ Ignace Aristide Minlend}
\address{\ African Institute for Mathematical Sciences (A.I.M.S.) of Senegal KM 2, Route de Joal, B.P. 1418 Mbour, S\'en\'egal }
\email{mouhamed.m.fall@aims-senegal.org, $\qquad $ignace.a.minlend@aims-senegal.org.}

\maketitle

\textbf{Abstract:} Let $(\M,g)$ be a compact Riemannian manifold of
dimension $N$, $N\geq 2$.  In this paper, we prove  that there exists
a family of domains $(\O_\e)_{\e\in(0,\e_0)}$   and functions $u_\e$
such that
\begin{align}\label{eq:absSerrin}
  \begin{cases}
    -\Delta_{g} u_\e=1& \quad \textrm{ in }  \Omega_\e\vspace{3mm}\\
u_\e=0&  \quad\textrm{ on }\partial\Omega_\e\vspace{3mm}\\
 {g}(\nabla_{ {g}} {u_\e}, {\nu}_\e)=-\frac{\e}{N}&   \quad \textrm{ on }\partial\Omega_\e,
  \end{cases}
  \end{align}
where   $\nu_\e$ is the unit outer normal of $\partial\Omega_\e$.
The domains $\O_\e$ are smooth perturbations of geodesic balls of radius $\e$. If, in addition,
$p_0$ is a non-degenerate critical point of the scalar curvature of
$g$ then, the family $(\de\O_\e)_{\e\in(0,\e_0)}$ constitutes a
smooth foliation of a neighborhood of $p_0$. By considering a family of  domains $\O_\e$  in which \eqref{eq:absSerrin}  is satisfied, we also prove that if this family
converges to some point $p_0$ in a suitable sense as $\e\to 0$, then $p_0$ is a critical point of the scalar curvature.
A Taylor expansion of he energy rigidity for the torsion problem is also given.\\

\bigskip
  \noindent 2010 {\it Mathematics Subject Classification.}   	58J05, 58J32,58J37, 35N10,   	35N25 .\\
  \noindent {\it Keywords.}  Over-determined problem, foliation, isochoric profile.
\section{Introduction }

    Let $(\M,g)$ be a compact Riemannian manifold of dimension $ N\geq2$.
    We are interested in this paper in the  construction of  smooth domains $\O\subset\M$ where there   exists  $u\in C^2(\ov{\O})$    such that
\begin{align}\label{eq:problemb}
  \begin{cases}
    -\Delta_{g} u=1& \quad \textrm{ in }  \Omega\vspace{3mm}\\
u=0&  \quad\textrm{ on }\partial\Omega\vspace{3mm}\\
{g}(\nabla_{ {g}} {u }, {\nu})=-c.&   \quad \textrm{ on }\partial\Omega,
  \end{cases}
  \end{align}
where $c$ is a positive constant,
$\Delta_{g}=\textrm{div}(\nabla_{g})$ is the Laplace-Beltrami
operator, and  $\nu$ is the unit outer  normal of $\partial\Omega$.
In the Euclidean space (at least in $\R^2$), if a function $v$
satisfies the first two equations of \eqref{eq:problemb} then the quantity
$$
\int_{\O}|\n v|^2
$$
is the torsional rigidity of the rod $\O\times \R$. Namely the torque required
for unit angle of twist per unit length.  We refer to \cite{Landau} for the precise derivation and its relation with bending a plane membrane and the motion of viscous fluids. Still in Euclidean space, it was  proved by Serrin in \cite{Serrin}
that a $C^2$ domain $\O$ in which \eqref{eq:problemb} has a solution
must be a ball. The argument of Serrin to
 prove his result relies on the moving plane method due to Alexandrov in \cite{Alexandrov}.
 In fact Alexandrov introduced the moving plane method while proving that an embedded constant mean curvature hyper-surface
 in $\R^N$  must be a sphere.   Serrin's result can be also derived from the Alexandrov's. Namely if \eqref{eq:problemb}
  has a solution then $\de\Omega$ has constant mean curvatures, see the work of Farina-Kawohl in \cite{FK} and   Choulli-Henrot\cite{CH}.

While CMC hyper-manifolds  are stationary sets for the area functional under volume
 preserving deformations, an over-determined problems arises when looking for a  stationary
 set (under volume preserving
 deformations as well) to some energy functional given by  some functional inequalities. In our case  this energy is   proportional to the inverse of the  torsional rigidity:
\begin{equation}\label{eq:propp99}
J(\Omega):=\inf\biggl\{\int_{\Omega}|\n
u|^{2}_g\,\textrm{dvol}_g:~~\int_{\Omega}u\,\textrm{dvol}_g=1,~u\in
H^{1}_{0}(\Omega)\biggl\}.
\end{equation}
 In particular minimizing $\O\mapsto J(\O)$ is equivalent to maximizing the torsion rigidity and therefore Serrin's result states that balls maximize the torsion rigidity as it can be also derived from the Faber-Krahn inequality, see for instance \cite{Chavel-1}.\\
A   smooth bounded domain $\O$ is stationary (or a critical point) for the functional $\O\mapsto J(\O)$ under volume preserving deformations if and only if there exists $u_\O\in H^1_0(\O)$ such that
\begin{equation}\label{eq:overdet0}
\left\{ \begin{array}{ll}
\displaystyle -\Delta_g u_{\Omega}=J(\Omega)&\quad \textrm{ in }~\Omega\vspace{3mm}\\
\displaystyle u_{\Omega}=0&\quad \textrm{ on }~ \partial \Omega\vspace{3mm}\\
\displaystyle  {g}(\nabla_{ {g}} {u_\O}, {\nu})=Const.&\quad \textrm{ on }~ \partial \Omega.
\end{array} \right.
\end{equation}
We refer to Section 3 for more detailed explanations. In Euclidean
space it is known from Serrin \cite{Serrin} and Weinberger's \cite{Wein} work that stationary
smooth domains are balls. In this paper, we will show that in a
Riemannian manifold, geodesic balls can be perturbed to stationary
sets for $J$. Before stating our result, let us recall some known
results in the constructions of  CMC hyper-manifolds. In \cite{Ye},  Ye proved   that if $p_0$ is a non-degenerate
critical point of the scalar curvature of $g$ then the geodesic ball
$B^g_\e(p_0)$ might be perturbed  to a CMC sub-manifolds with mean
curvature equal $\frac{N-1}{\e}$, the mean curvature of Euclidean balls $B_\e(0)$ with radius $\e$.
   By taking advantages on the variational properties of the problem,   Pacard and Xu showed
in \cite{PX} that CMC hyper-manifolds with mean curvature $ \frac{N-1}{\e}$ always
exist and the number is not less than the category of $\M$.

 By direct computation, a solution to \eqref{eq:problemb} in $\R^N$ is given by
$$\phi_0^\e(x):=\frac{\e^2-|x|^2}{2N}$$ which clearly satisfies
 \begin{align}\label{eq:ropp2I}
  \begin{cases}
    -\Delta\phi_0^\e=1& \quad \textrm{ in }  B_\e(0) \vspace{3mm}\\
\phi_0^\e=0&  \quad\textrm{ on }\partial B_\e(0) \vspace{3mm}\\
\dfrac{\partial \phi_0^\e}{\partial\nu }=-\frac{\e}{N}&   \quad \textrm{
on }\partial B_\e(0).
  \end{cases}
  \end{align}

The
main results in this paper parallel those of Ye and Pacard-Xu. We first prove the following:

\begin{Theorem}\label{theo0}
 Let $(\M,g)$ be a compact Riemannian manifold of dimension $N\geq2$. There exists $\e_{0}>0$ and  a smooth function
$\calF: \M\times(0,\e_0)\longrightarrow \mathbb {R}$ such that
 for all $\e\in (0,\e_0)$, if $p$
is a critical point of the function $\calF(.,\e)$ then  there
exists a smooth domain $\Omega_{\e}$ and a   function
$u_{\e}\in C^2(\ov{\O_\e})$  such that
\begin{align}\label{eq:problemz}
  \begin{cases}
    -\Delta_{g} u_\e=1& \quad \textrm{ in }  \Omega_\e \vspace{3mm}\\
u_\e=0&  \quad\textrm{ on }\partial\Omega_\e \vspace{3mm}\\
{g}(\nabla_{ {g}} {u_\e}, {\nu}_\e)=-\frac{\e}{N}&   \quad \textrm{
on }\partial\Omega_\e.
  \end{cases}
  \end{align}
 Here  $\nu_\e$ stands for  the unit outer normal of $\partial\Omega_\e$.
  Moreover we have
\be \label{eq:Phiep}
 ||\calF(\cdot,\e)-S_g ||_{C^{2,\a}(\M)}\leq C\e^{2},
\ee
where $S_g$ is the scalar curvature of $\M$ and  $C>0$ is a constant independent on $\e$.
  \end{Theorem}
Let us denote by $B_\e^g(p)$ the geodesic ball centered at $p$ with
radius $\e$. The domains $\Omega_\e$ we construct are perturbations
of geodesic balls in the sense that $ \Omega_{\e}
=(1+v^{p,\e}) B_\e^g(p)$,  with $v^{p,\e}:\de B_\e^g(p) \to \R $
satisfying
$$||v^{p,\e}||_{C^{2,\alpha}(\partial
B^{g}_{\e}(p))}\leq c\e^{2}$$
while the function $u_\e$ satisfies the estimates
\be
\|u_\e\|_{C^{2}(\ov{ \Omega_{\e}})}\leq c,
\ee
where the constant $c$ is independent on $\e$.\\
The fact that a solution to \eqref{eq:problemz} exists is guaranteed
by the existence of $\calF$. Indeed the number   of critical
points of $p\mapsto \calF(p,\e)$ is greater than the
Lusternik-Shnirelman category of $\M$,  see \cite{AM}. However
\eqref{eq:Phiep}  implies that near a  topologically stable critical
point of $S_g$, there exists a critical point of $\calF(\cdot,\e)$.
In particular  if $p$ is  a non-degenerate critical point of $S_g$
then the implicit function theorem implies that there exists a curve
$p_\e$ of critical points of $\calF(\cdot,\e)$. It is known from the
work of Micheletti and Pistoia in \cite{MiPi}  that for a generic
metric on a manifold, all critical points of the scalar curvature
are non-degenerate.  This implies that  for a generic metric $g'$, a
neighborhood of  any critical point of $S_{g'}$ can be foliated by
CMC hyper-manifolds, similar to  geodesic spheres,  thanks to Ye's
result. The analogous to this result is contained in the following
 \begin{Theorem}\label{theo1}
Assume that $p_0$ is a non-degenerate critical point of the scalar
curvature function $S_g$ of $(\M,g)$. Then, there exists
$\e_{0}>0$ such that $(\de\O_\e)_{\e\in(0,\e_0)}$ constitutes
a smooth foliation of a neighborhood of $p_0$, where $\O_\e$ is a domain in which Serrin's over-determined problem \eqref{eq:problemz}
possesses a solution.
\end{Theorem}
In fact we obtain a precise form of the boundary of the domains constructed  in Theorem \ref{theo1}. Indeed, we proved that  there exists a   function $\o^{\e}: S^{N-1}\to \R_+$ such hat
$$
\de \O_\e=\left\{\exp_{p_0}\left(\o^{\e}(y)\sum_{i=1}^N y^i E_i\right),\quad y\in S^{N-1} \right\}
$$
and moreover the map $\e\mapsto \o^\e$ satisfies $\de_\e\o^\e|_{\e=0}=1$.  In particular, we can see that  the domains   $\O_\e\subset B^g_{\d_\e}(p_0)$, for some function $\d_\e=\e+O(\e^2)$.
Our next result can be merely seen as    the converse of Theorem \ref{theo1}.

\begin{Theorem} \label{theo2}
Suppose that  for every $\e>0$ there exist  $\d_\e >0$, a smooth domain $\O_\e\subset B^g_{\d_\e}(p_0)$ and a function $u_{\e}\in C^2(\ov{\O_\e})$  such that
\begin{align}\label{eq:problemz00}
  \begin{cases}
    -\Delta_{g} u_\e=1& \quad \textrm{ in }  \Omega_\e \vspace{3mm}\\
u_\e=0&  \quad\textrm{ on }\partial\Omega_\e \vspace{3mm}\\
{g}(\nabla_{ {g}} {u_\e}, {\nu}_\e)=-\frac{\e}{N}&   \quad \textrm{
on }\partial\Omega_\e.
  \end{cases}
  \end{align}
Assume that
\be\label{eq:ede}
\e^{-1} |\d_\e-\e|\to 0 \qquad \textrm{ as } \e \to 0
\ee
and
\be \label{eq:uC2bound}
\|u_\e\|_{C^{2}(\ov{ \Omega_{\e}})}\leq C,
\ee
for some $C>0$ independent on $\e$. Then $p_0$ is a critical point of the scalar curvature $S_g$.
\end{Theorem}


An other question of interest we study in this paper is the expansion of the isochoric profile corresponding to the torsion problem. We define the profile  $\calT_\M$   by
$$
\calT_\M(v,g):= \inf_{|\O|_g=v} J(\O).
$$
In particular,  thanks to the Faber-Krahn inequality,
$$
\calT_{\R^N}(v)=J(B_1)\biggl(\frac{|B_1|}{v}\biggl)^{-\frac{N+2}{N}}.
$$
\begin{Theorem}\label{th:IM}
 We have
$$
 \calT_\M(v,g)=\left( 1- c_N\, v^{\frac{2}{N}}\,\max_{\M}S_g +O(v^{\frac{3}{N}} )  \right)
 \calT_{\R^N}(v),
$$
as $v\to 0$, where
$$
 c_N= \frac{N+6}{6N(N+4)}|B_1|^{-\frac{2}{N}}.
 $$
 \end{Theorem}
 This result suggests that torsion rigidity is maximized by sets
located where  scalar curvature is maximal. Let  $g_k$ be a metric of constant sectional curvature $k$ on a manifold $\M'$ with dimension $N$. Suppose that $ \max_{\M}S_g<N(N-1)k$ then Theorem \ref{th:IM} implies that
$$
 \calT_\M(v,g)> \calT_{\M'}(v,g_k)= J(B^{g_k}_v),\qquad \textrm{ as } v\to 0,
$$
where $ B^{g_k}_v$ is a geodesic ball ball with volume $v$ in $(\M, g_k)$.
We  quote    \cite{GLL} and \cite{Xiao} for some recent geometric comparisons of the energy torsional rigidity.

\begin{Remark}
The result in this paper provides critical domains that concentrate
at points. In a forthcoming work, we will be interested in
concentrations at minimal submanifolds. Namely
  letting $K$ be a \textit{non-degenerate minimal submanifold} of $\M$. Let $K_\rho$ be the geodesic neighborhood of $K$ with radius $\rho>0$.
We will perturbe the tube $K_\rho$   to a domain $\O_\rho$  such
that  there exists a function $u_\rho$  which satisfies
\begin{align*}
\begin{cases}
\displaystyle-\D u_\rho=1&\quad\textrm{ in } \O_\rho \vspace{3mm}\\
\displaystyle u_\rho=0&\quad\textrm{ on } \de\O_\rho\vspace{3mm} \\
%
%
 {g}(\nabla_{ {g}} {u_\rho}, {\nu}_\rho)=Const.& \quad\textrm{ on } \de \O_\rho.
\end{cases}
\end{align*}
In the CMC theory,  minimal submanifolds play as well an important
role.
 In comparison to Ye's result, nondegenracy of critical point of the scalar curvature is replaced by the fact $K$
 is non-degenerate: the Jacobi operator about $K$ does not have zero eigenvalues. We might not expect existence
 for all $\rho>0$ small but a sequence of $\rho's$ will do.  This is motivated by the work of  Malchiodi and Montenegro in \cite{malm}
  and related works on CMC's concentrating along submanifolds, \cite{mmm,mp,mmp}.
\end{Remark}

 Construction of solutions to
over-determined problems on Riemannian manifolds was first studied by Pacard and Sicbaldi in \cite{FPP}, where
they study an over-determined problem for the first Dirichlet
eignevalue $\l_1(\O_\e)$:
\begin{align}\label{eq:problemlam}
  \begin{cases}
    -\Delta_{g} u_\e=\l_1(\O_\e) \,u_\e & \quad \textrm{ in }  \Omega_\e\\
u_\e=0&  \quad\textrm{ on }\partial\Omega_\e\\
 {g}(\nabla_{ {g}} {u_\e}, {\nu}_\e)=Const.&   \quad \textrm{ on }\partial\Omega_\e.
  \end{cases}
  \end{align}
Pacard and Sicbaldi in \cite{FPP} proved  that  when the Riemannian
manifolds has a non-degenerate critical point $p_0$ of the scalar
curvature then it is possible to build extremal domains for any
given volume small enough, and such domains are close to geodesic
balls centered at $p_0$. This result has been improved by Delay and Sicbaldi
\cite{DS} eliminating the hypothesis of non-degenerate critical
point of the scalar curvature. In particular they showed the
existence of extremal domain of small volume in any
compact Riemannian manifold. Some other results and works on construction of solutions to  over-determined problems on Riemannian manifolds can be found
in \cite{DS,MS,ScSi,dPPW,Sic}.\\

We shall now explain our argument of proof which is based on geometric variational arguments,  see   the work of \cite{Kap},  \cite{PX, mmf,mmfcm1,mmfcm}, for the construction of
constant mean curvature hyper-surfaces and Delay-Sicbaldi \cite{DS}, for
the construction of extremal domains for the first eigenvalue of the
Laplace-Beltrami operator. See also \cite{AB} and \cite{AM} for related abstract perturbative methods. \\

 The idea is to perturbed a geodesic ball $B^g_\e(p)$. For any function $v\in L^2(S^{N-1})$, we will consider the decomposition $v=v_0+\bar v$ where $  \int_{S^{N-1}}\bar v \,\textrm{dvol}_{S^{N-1}}=0$.
 We define the scaled metric $\bar g= \e^{-2}g$. For $(v_0,\bar v)\in\R\times C^{2,\a}(S^{N-1})$ we consider the nearby sets of  the rescaled ball $B^{\bar g}_\e(p)$ given by:
 $$
B_{1+v}^{\bar g}(p):= \left\{ \exp^{ \bar{g}}_{p}
\biggl(\biggl(1+v_{0}+\chi(x)\bar{v}(x/|x|)\biggl)
\sum_{i=1}^Nx^{i}E_{i}\biggl)\,:\, |x|<1  \right\},
 $$
 where $\chi$ is a radial  cutoff function identically
equal to $0$ for $|x|\leq\frac{1}{4}$ and $1$ for
$|x|\geq\frac{1}{2}$.\\

The main idea is to find $p, v_0,\bar v$ such that Serrin's over-determined problem \eqref{eq:problemb} is solvable.
The first step consists in constructing a    first approximate solution by solving only the Dirichlet problem:
given a point $p\in \M$, there exists $\e_0>0$ such that for
all $\e \in(0,\e_0)$  and for all $(v_0,\bar{v})\in
\mathbb{R}\times C^{2,\alpha}(S^{N-1})$ satisfying
$$
|v_0|\leq\e_0,\qquad\qquad  ~~||\bar{v}||_{C^{2,\alpha}(S^{n-1})}\leq\e_0~~\qquad \textrm{ and } \qquad ~~\int_{S^{N-1}}\bar{v}\, \textrm{dvol}_{S^{N-1}}=0,
$$
 there exists a unique positive function $\bar{\phi}=\bar{\phi}(p,\e,v_0,\bar{v}) \in C^{2,\alpha}(B^{\bar{g}}_{1+v}(p))$ such that
\begin{equation}\label{eq:prop1I}
\left\{ \begin{array}{ll} -\Delta_{\bar{g}} \bar{\phi}=1& ~~\quad \textrm{ in }
~B^{\bar{g}}_{1+v}(p) \vspace{3mm}\\
\bar{\phi}=0& ~~\quad \textrm{ on } ~~\partial B^{\bar{g}}_{1+v}(p),
\end{array} \right.
\end{equation}
This is done removing the dependence of the domains on the
parameters by considering a change of variable via the  function
$Y_{p,v} $
$$
Y_{p,v}(x):=  \exp^{ \bar{g}}_{p}
\biggl(\biggl(1+v_{0}+\chi(x)\bar{v}(x/|x|)\biggl)
\sum_{i=1}^N x^{i}E_{i}\biggl)
$$
which parameterizes $B_{1+v}^{\bar g}(p) $ over the unit  ball $B_1$ centered at the origin. Hence with the pull-back metric  $\hat g$ of $\bar g$ with respect to $Y_{\e,v}$, \eqref{eq:prop1I} becomes
\begin{equation}\label{eq:prop1Ih}
\left\{ \begin{array}{ll} -\Delta_{\hat{g}} \hat{\phi}=1& ~~\quad \textrm{ in }
~B_1 \vspace{3mm}\\
\hat{\phi}=0& ~~\quad \textrm{ on } ~~\partial B_1.
\end{array} \right.
\end{equation}
Once we find $v_0,\bar v$ so that  \eqref{eq:prop1Ih} holds for all $p\in\M$, we compute the normal derivative of $\hat {\phi}$.
Denote by $\hat  \nu$ the unit outward normal to $B_1$ with respect to the metric $\hat g$.
We have  obtained
\begin{align*}
\hat{g}(\nabla_{\hat{g}}\hat{\phi},\hat{\nu})_{\mid_{\partial
B_1}}=-|\n_g\hat{\phi} |_{\hat{g}}=-\frac{1}{N}+\frac{1}{N}\left[(\partial_{\nu}\psi_{v})_{|_{\partial
B_1}}- {v}\right]  + error,
\end{align*}
where $\psi_v$ satisfies
\begin{equation}\label{eq:pemvwI}
\displaystyle\left\{ \begin{array}{ll} \Delta \psi_{\bar v}=0& \quad~~\textrm{ in }~~B_1\vspace{3mm} \\
\displaystyle\psi_{\bar v}= \bar{v}&\quad ~~\textrm{ on }
~~\partial B_1.
\end{array} \right.
\end{equation}
The second step is to find $p,v_0,\bar v$ such that
\be\label{eq:Geq0}
G(\e,p,v_0,\bar v):=\frac{1}{N}\left[(\partial_{\nu}\psi_{v})_{|_{\partial
B_1}}- {v}\right]+error=0.
\ee
Direct computations then give
$$
\frac{\de G}{\de (v_0,\bar{v})}(0,p,0,0)[w_0,\bar w]= \frac{1}{N}\left[(\partial_{\nu}\psi_{w})_{|_{\partial
B_1}}- {w}\right].
$$
Note that  the map $\bar v\mapsto \partial_{\nu}{\psi_{\bar v}}_{|_{\partial
B_1}}$ is the classical Dirichlet-to-Neuman operator.
Its spectrum is known and the eigenvalues are the spherical harmonics.  It is then easy to verify that
$$
Ker \frac{\de G}{\de (v_0,\bar{v})}(0,p,0,0)=\{ x^i\,:\, i=1,\dots,N
\}.
$$
This allows us to solve \eqref{eq:Geq0} modulo its kernel. Namely
there exist  $v^{\e,p}=v_0^{\e,p}+\la a^{\e,p},x\ra+ \bar{v}^{\e,p}$  such that
\be\label{eq:Geq01} G(\e,p,v_0^{\e,p},\bar v)=- \la a^{\e,p},x
\ra,\quad \forall x\textrm{ in } S^{N-1} . \ee
Gathering what we have so far, we may say that we have found a
function $\bar{\phi}^{\e,p}$ satisfying
\begin{equation}\label{eq:problem-aI}
\left\{ \begin{array}{ll} -\Delta_{\bar g} \bar{\phi}^{\e,p}=1& ~~\textrm{ in } ~~B^{\bar{g}}_{1+v^{\e,p}}(p)\vspace{3mm}\\
\bar{\phi}^{\e,p}=0 &~~\textrm{ on }  ~~\partial B^{\bar{g}}_{1+v^{\e,p}}(p)\vspace{3mm}\\
\bar{g}(\nabla_{\bar{g}}\bar{\phi}^{\e,p},\bar{\nu})=-\frac{1}{N}-\bar{g}(A^{\e,p},{\calV}^{\e,p})&~~\textrm{on}
~~\partial B^{\bar{g}}_{1+v^{\e,p}}(p),
\end{array} \right.
\end{equation}
where for all $x\in S^{N-1}$, we define $
A^{\e,p}(Y_{p,v^{\e,p}}(x)):=dY_{p,v^{\e,p}}(x)[ a^{\e,p}] $ and
similarly    ${\calV}^{\e,p}:=dY_{p,v^{\e,p}}(x)[ x]  $. This program is detailed in Paragraph \ref{ss:fixedpoint}.

Let us remind that the domains we are looking for are critical
points of  the energy functional $J(\O)$ under volume constraints
and thus  by the Lagrange multiplier rule, they are critical points
of $J(\O)+\l |\O|_{\bar{g}}$, for some $\l\in\R$. We will take this
advantage in the  third step to annihilate $A^{\e,p}$ for some
special points $p$. Indeed by defining
$$
\Phi_\e(p)=J( B^{\bar{g}}_{1+v^{\e,p}}(p))+\frac{1}{N^2}
|B^{\bar{g}}_{1+v^{\e,p}}(p)|_{\bar g},
$$
we will show that if $p$ is a   critical points of this functional
then $A^{\e,p}=0 $. Rescaling back, we get  the desired result:
there exists a  function $u_\e $ such that
\begin{equation}\label{eq:problem-aI}
\left\{ \begin{array}{ll} -\Delta_{  g} u_\e=1& ~~\textrm{in} ~~B^{ {g}}_{\e(1+v^{\e,p})}(p)\vspace{3mm}\\
u^\e=0 &~~\textrm{on}  ~~\partial B^{ {g}}_{\e(1+v^{\e,p})}(p)\vspace{3mm}\\
 {g}(\nabla_{ {g}}u^\e, {\nu}_\e)=-\frac{\e}{N} &~~\textrm{on}
~~\partial B^{ {g}}_{\e(1+v^{\e,p})}(p).
\end{array} \right.
\end{equation}
We refer to  Paragraph \ref{ss:Geometric-var} for more details.
In addition the functional  $\Phi_\e$   has a Taylor expansion for which the main term is given by the scalar curvature, see Lemma \ref{Lem4.3}.\\

 Next, in Section \ref{s:local-fol}, we will prove that   we have a smooth foliation near non-degenerate critical points of the scalar curvature.
 Here  we take advantages of the expansion of $\Phi_\e(p)=\a_n+\b_n \e^2S_g(p)+O(\e ^4) $ to see that provided $p_0$ is
  a non-degenerate critical point of the scalar curvature   there exists a curve $p_\e$ of critical points of $\Phi_\e$ such that $\textrm{dist}_g(p_\e,p_0)\leq C \e^2$.
This fact allows us to re-parameterize $\partial B^{
{g}}_{\e(1+v^{\e,p})}(p)$ by perturbed sphere with increasing radius
$\o^\e$. Indeed there exists a   nonnegative function $\o^\e$
such that
$$
\partial B^{ {g}}_{\e(1+v^{\e,p_\e})}(p_\e)=\left\{\exp_{p_0}\left( {\o}^\e(y)  y^i E_i\right)\,:\, y\in S^{N-1} \right\},
$$
with
$$
\o^\e(\cdot)>0,\qquad \textrm{ and } \qquad  \de_{\e}
{\o}^\e (\cdot)>0.
 $$
 It is worth noticing that from our argument to
prove local foliation, the sets $\O_\e$ in \eqref{eq:problemlam}
constructed by Pacard and Sicbaldi in \cite{FPP} enjoys such a local
foliation, see Remark \ref{rem:foliation-PS}.\\

Finally in Section \ref{s:converse}, we prove  Theorem \ref{theo2}. The proof is based on  the regularity result of nearly minimizing sets for the perimeter functional. Indeed, just integrating the \eqref{eq:problemz00}, we see that the domains $\O_\e$ satisfies
$$
|\de\O_\e|_g=\frac{N}{\e} | \O_\e|_g
$$
while  \eqref{eq:ede} shows that they are contained in the  ball  $(\e+o(\e)) B_1$. This implies that
$$
  |\de \O_\e|_g\leq
 (1+o(1))c_N |\O_\e|_g^\frac{N-1}{N},
$$
where $c_N=N|B|^\frac{1}{N}$ is the isoperimetric constant of  ${\R^{N}}$. Therefore up to a scaling they   nearly minimize the area functional among domains with   volume $|B_1|$. Using some simple arguments, we deduce that they have bounded boundary mean curvatures. This leads to smooth convergence to the unit ball.  We note that even if our argument works also when considering CMCs instead of critical domains,  we choose not to expose it here.  Among others we quote \cite{OD}, \cite{Nard}, \cite{Sun}, \cite{Laurin} and  \cite{mp}, where the authors   characterized the sets where a sequence of CMC's hyper-surfaces converges as their mean curvature tends to infinity. From the work of \cite{Sun} and  \cite{Laurin},  it is also naturel to expect that the assumption \eqref{eq:ede} can be relaxed.

\bigskip
\noindent
\textbf{ Acknowledgements}: This work is supported by the Alexander von Humboldt foundation and the    German Academic
Exchange Service (DAAD).

\section{Preliminaries and notations }
Given a point  $p\in \M$, we let   $E_{1},...,E_{N}$ be an
orthonormal basis of the tangent plane
 $T_{p}\M$. We consider geodesic coordinate
system
$$
\mathbb{R}^{n} \ni (x^1,...,x^N)=x \longmapsto
F_p(x):=\exp_{p}(X)\in \M,
$$
 where we use here and in the following  the notation
  $$
X:=\sum^{N}_{i=1}x^{i}E_i\in T_{p}\M.
$$
 The map $F_p$ induces
coordinate vector fields
 $$X_i:= dF_p (x)[E_i].$$
 We denote by
$$
R_p:T_p\M\times T_p\M\times T_p\M\longrightarrow T_p\M
$$
the Riemanniann curvature tensor at $p$ and
 $$
 Ric_p:T_p\M\times
T_p\M\longrightarrow
\mathbb{R}, \qquad  Ric_p(X,Y)=-\sum^{N}_{i=1}g\biggl(R_p(X,E_i)Y,E_i)\biggl)
$$
 the Ricci curvature tensor at $p$.
The scalar  curvature of $(\M,g)$ at $p$ is defined  by
$$
S_g(p)=\sum^{N}_{k=1}Ric_p(E_k,E_k).
$$
At a point $q=\exp_p(X) $, we  define
$$g_{ij}(x):=g(X_i,X_j).$$
The proof of the expansion of the metric $g$ near $p$ in normal
coordinates is classical and can be found in   \cite{TJW}
or \cite{STY}.
\begin{Proposition} \label{prop:expgij}
At a point  $q=\exp_p(X)$, we have
$$g_{ij}(x)=\delta_{ij}+\frac{1}{3} g(R_p(E_i,X,)E_j,X) +\frac{1}{6} g(\n_XR_p(E_i,X) E_j,X) +O_{p}(|x|^{4}),$$
as $|x|\to 0$.
\end{Proposition}
%
%
Let $f: S^{N-1}\longrightarrow (0, \infty)$ be a continuous function
whose $L^{\infty}$ norm is small (say less than the cut locus of
$p$). We can decompose $f$  into $f=f_0+\bar{f}$, where $f_0$ is a
constant  and $\bar{f}$ has mean value equal to $0$. We define
$$
B^{g}_{f}(p):=\biggl\{\exp^{g}_{p}\left((f_0+\chi\bar{f}(x/|x|))X\right)\quad :\quad
|x|<1\biggl\},
$$
 where $\chi$ is a radial  cutoff function identically
equal to $0$ for $|x|\leq\frac{1}{4}$ and $1$ for
$|x|\geq\frac{1}{2}$. In particular if $f_0\equiv r$ a positive constant and $\bar f=0$, then
$ B^g_f(p)$ is nothing but the geodesic ball centered at $p$ with
radius $r$.

Similarly, we denote by $|\O|_g$ the volume in the metric $g$ of a
smooth domain $\O\subset\M$, $\textrm{dvol}_g$ the volume element in
the metric $g$ to integrate over the domain and $d\sigma
_g$ denotes
the volume element in the induce metric $g$ to integrate over the
boundary of a domain. $\Delta_{g}$ and $\nabla_{g}$ denotes
respectively, the Laplace-Beltrami and the gradient operator with
respect to the metric $g$. It will be understood that when we do not
indicate the metric as a
superscript, we understand that we are using the Euclidian one.\\

    Our aim is to show that, for $\e>0$ small enough, we can find
a point $p$ and a (small) function $v:S^{N-1}\longrightarrow (0,
\infty)$ such that, on $(\calM,g),$ the over-determined problem

\begin{align}\label{eq:problem00}
 \begin{cases}
 \displaystyle -\Delta_{g} u=1& ~~\textrm{in}~~B^{g}_{\e(1+v)}(p)\vspace{3mm}\\
 \displaystyle u=0& ~~\textrm{on} ~~\partial B^{g}_{\e(1+v)}(p)\vspace{3mm}\\
  \displaystyle g(\nabla_{g}u,\nu)=-\frac{\e}{N}&~~\textrm{on} ~~\partial B^{g}_{\e(1+v)}(p)
\end{cases}
\end{align}
has a solution, where
$\nu$ is the unit outer normal vector about $\partial B^{g}_{\e(1+v)}(p)$.\\

We consider the dilated
metric $\bar{g}=\e^{-2}g$ and   rewrite \eqref{eq:problem00} on $(\calM,\bar g),$  as
$$
\left\{ \begin{array}{ll} -\Delta_{g} \bar{u}=1& ~~\textrm{in} ~~B^{\bar{g}}_{1+v}(p)\vspace{3mm}\\
\bar{u}=0 &~~\textrm{on}  ~~\partial B^{\bar{g}}_{1+v}(p)\vspace{3mm}\\
\bar{g}(\nabla_{\bar{g}}\bar{u},\bar{\nu})=-\frac{1}{N}&~~\textrm{on} ~~\partial
B^{\bar{g}}_{1+v}(p),
\end{array} \right.
$$
where
$$\bar{u}=\e^{-2}u .$$
The Taylor expansion of the scaled metric $\bar{g}$ can be easily derived from Proposition \ref{prop:expgij}. Indeed we have
\begin{equation}\label{exp}
\bar{g}_{ij}(x)=g_{ij}(\e x)=\delta_{ij}+\frac{\e^2}{3}
g(R_p(E_i,X)E_j,X) +\frac{\e^3}{6} g(\n_XR_p(E_i,X) E_j,X)
+O_{p}(\e^{4}).
\end{equation}
Given $v\in C^{2,\a}(S^{N-1})$, with $\a\in (0,1)$, we can decompose $v$ as $v=v_0+\bar v$, where
$$
\int_{S^{N-1}}\bar{v}\,dvol_{S^{N-1}}=0.
$$
The perturbed geodesic ball $ B^{\bar{g}}_{1+v}(p)$ can be
parameterized by the map $Y_{p,v}: B_1\to  B^{\bar{g}}_{1+v}(p)$
given by
\be \label{eq:Yv} Y_{p,v}(x):=\exp^{ \bar{g}}_{p}
\biggl(\biggl(1+v_{0}+\chi(x)\bar{v}(x/|x|)\biggl)\sum^{N}_{i=1}
x^{i}E_{i}\biggl).
\ee
In the following, we will put   $\rho:=1+v$ and  denote by
$\rho_i$ (resp. $\rho_{ij}$) the partial derivative of $\rho$ with
respect to $x^i$ (resp.  the partial derivatives with respect to $x^i$ and $x^j$).\\
The parametrization \eqref{eq:Yv} induces a metric $\hat{g}$ on $\R^N$.
 Our next task is to derive the Taylor expansion of the metric $\hat g $.
 To this end, we will need to fix some notations.\\

\noindent \textbf{Notations:} Any expression of the form $L_p^i(v)$
denotes a linear combination of the function $v$ together with its
partial derivatives with respect to $x^i$ up to order $i=0,1,2$. The
coefficient of $L_p^i$ might depend on $\e$ and $p$ but, for all
$k\in\mathbb{N}$, there exists a constant $c>0$ independent of
$\e$ and $p$ such that
$$||L_p^i(v)||_{C^{k,\alpha}(S^{N-1})}\leq
c||v||_{C^{k+i,\alpha}(S^{N-1})}.$$ Similarly,   any expression of the form $Q^i_p(v)$ denotes a
nonlinear operator in the function $v$ together with its derivatives
with respect to  $x^i$ up to order $i$. The coefficient of the
Taylor expansion of $Q^i_p(v)$ in power of $v$ and its partial
derivatives might depend on $\e$ and $p$ and, given $k\in
\mathbb{N}$, there exists a constant $c>0$ such that
$Q^i_p(0)=0$ and
\begin{align*}
||Q^i_p(v_1)-Q^i_p(v_2)||_{C^{k,\alpha}(S^{N-1})}&\leq
c\biggl(||v_1||_{C^{k+i,\alpha}(S^{N-1})}+||v_1||_{C^{k+i,\alpha}(S^{N-1})}\biggl) \times\\
||v_1-v_2||_{C^{k+i,\alpha}(S^{N-1})},
 \end{align*} provided
$||v_1||_{C^{1,\a}(S^{N-1})}+||v_2||_{C^{1,\a}(S^{N-1})}\leq 1$. Terms of the form
$O_{p}(\e^{l})$ are smooth functions on $S^{N-1}$ that might
depend on $p$ but which are bounded by a constant (independent of
$p$) times $\e^{l}$ in the $C^{k}$ topology,
for all $k\in \mathbb{N}$. Finally the function $P_\e^i(v)$ stands for
$$
P_\e^i(v)=\e^2 L^i(v)+ Q^i(v)+O_p(\e^{4}).
$$

We recall that the map $Y_{p,v}$  parameterizes
$B^{\bar{g}}_{1+v}(p)$ and we denote  by $\hat{g}$ the pull-back metric on $B_1$ via
 $Y_{p,v}$. At
the point $q=Y_{p,v}(x)$, we define
$$\hat{g}_{ij}(x):=\bar{g}\left(\frac{\de Y_{p,v}}{\de
x^i}(x),\frac{\de Y_{p,v}}{\de x^j}(x) \right) .$$

\begin{Lemma}\label{Lem2} For all $x\in B_1$, we have the following expansions
$$\hat{g}_{ij}(x)=\rho^{2}\biggl(\delta_{ij}+\rho_ix^j+\rho_jx^i+\frac{\e^2}{3} g(R_p(E_i,X) E_j,X) +\frac{\e^3}{6} g(\n_XR_p(E_i,X) E_j,X)+P_{\e}^1(v)\biggl)$$
and
\begin{align*}
\Delta_{\hat{g}}&=\rho^{-2}\Delta-2\sum^{N}_{i,j=1}x^i\rho_{j}\partial^{2}_{ij}-2\sum^{N}_{j=1}\rho_{j}\partial_j-\Delta\rho
\sum^{N}_{j=1}x^j\partial_j-\frac{\e^{2}}{3} \sum^{N}_{i,j=1}g(R_p(E_i,X) E_j,X)\partial^{2}_{ij}\\
&\quad +\frac{2\e^{2}}{3} \sum^{N}_{i,j=1}g(R_p(E_i,E_j)
E_i,X)\partial_j+\frac{\e^{3}}{3}\sum^{N}_{i,j=1}g(\n_XR_p(E_i,X)
E_i,E_j)\partial_{j}\\
&\quad -\frac{\e^{3}}{6}\sum^{N}_{i,j=1}g(\n_{E_i}R_p(E_i,X)
E_j,X)\partial_{j}+\frac{\e^{3}}{12}\sum^{N}_{i,j=1}g(\n_{E_j}R_p(E_i,X)E_i,X)\partial_{j}\\
&\quad-\frac{\e^{3}}{6}\sum^{N}_{i,j=1}g(\n_XR_p(E_i,X)
E_j,X)\partial^{2}_{ij}+\sum^{N}_{i,j=1}\Delta_{\e,\bar{v}}^{ij},
\end{align*}
where
$$
 \Delta_{\e,\bar{v}}^{ij}=P_{\e}^2(v)\de_{ij}+P_{\e}^2(v)\partial_j.
$$
\end{Lemma}

\proof We have
$$\frac{\de Y_{p,v}}{\de x^i}(x)=\rho_{i}\sum^{N}_{k=1}x^{k}X_{k}+\rho X_{i}=\rho_{i}\Upsilon+\rho X_{i}~~\quad \forall i=1,...,
N,$$ ~~ where~~
\be\label{eq:Upsilon}
\Upsilon=\sum^{N}_{k=1}x^{k}X_{k}.
\ee
 We find using
the expansion \eqref{exp} that
$$
 \bar{g}(\Upsilon,
\Upsilon)\equiv |x|^{2} \qquad \textrm{ and }\qquad
 \bar{g}(\Upsilon,X_i)\equiv x^i,\quad i=1,\dots,N.
$$
These equalities then yield
\be \label{eq:ghatij}
\hat{g}_{ij} =\rho^{2}\biggl(\delta_{ij}+\rho_ix^j+\rho_jx^i+\frac{\e^2}{3} g(R_p(E_i,X, E_j,X) +\frac{\e^3}{6} g(\n_XR_p(E_i,X, E_j,X)+P_{\e}^1(v)\biggl).
\ee
 The  first expansion in the lemma then
follows.\\
The  expansion of Laplace-Beltrami operator of the metric $\hat{g}$
is obtained using the formula
\be\label{eq:formDel}
\Delta_{\hat{g}}=\hat{g}^{ij}\partial^{2}_{ij}+(\partial_{i}\hat{g}^{ij})\partial_{j}+\frac{1}{2}\hat{g}^{ij}(\partial_{i}\log|\hat{g}|)\partial_{j}=(1)+(2)+(3).
\ee We start with the last term. Thanks to \eqref{eq:ghatij}, it is
not difficult to see that
 \be \label{eq:ghat-inv-ij}
 \hat{g}^{ij}
=\rho^{-2}\biggl(\delta_{ij}-\rho_ix^j-\rho_jx^i-\frac{\e^2}{3}
g(R_p(E_i,X, E_j,X) -\frac{\e^3}{6} g(\n_XR_p(E_i,X,
E_j,X)+P_{\e}^1(v)\biggl).
 \ee
We also have
$$
\log|\hat{g}|=2N\log\rho+2\sum^{N}_{s=1}x^s\rho_s-\frac{\e^{2}}{3}Ric_{p}(X,X)+\frac{\e^{3}}{6}
\sum^{N}_{s=1}g(\n_XR_p(E_s,X)E_s,X)+P_{\e}^1(v)
$$
and  by a computation, we get
\begin{align*}
\partial_{i}(\log|\hat{g}|)&=2(N+1)\rho_i+2\sum^{N}_{s=1}x^s\rho_{is}+\frac{2\e^{2}}{3}\sum^{N}_{k=1}g(R_{p}(E_k,E_i)E_k, X)\\
&+\frac{\e^{3}}{3}\sum^{N}_{s=1}g(\n_XR_p(E_s,X)
E_s,E_i)+\frac{\e^{3}}{6}\sum^{N}_{s=1}g(\n_{E_i}R_p(E_s,X)
E_s,X)+P_{\e}^2(v).
\end{align*}
This together with \eqref{eq:ghat-inv-ij}   give
\begin{align*}
(3)&=(N+1)\sum^{N}_{j=1}\rho_j\partial_j+\sum^{N}_{ij=1}x^i\rho_{ij}\partial_j+\frac{\e^{2}}{3}\sum^{N}_{s,j=1}g(R_{p}(E_s,E_j)E_s, X)\partial_{j}\\
&+\frac{\e^{3}}{6}\sum^{N}_{s,j=1}g(\n_XR_p(E_s,X)E_s,E_j)\partial_{j}+\frac{\e}{12}\sum^{N}_{s,j=1}g(\n_{E_j}R_p(E_s,X)E_s,X)\partial_{j}+\sum^{N}_{j=1}P_{\e}^2(v)\partial_{j}.
\end{align*}
We compute the partial derivative of $\hat{g}^{ij}$ with respect to
$x_i$ and get
\begin{align*}
(2)&=-(N+3)\sum^{N}_{j=1}\rho_j\partial_j-\Delta\rho
\sum^{N}_{j=1}x^j\partial_j-\sum^{N}_{i,j=1}x^i\rho_{ij}\partial_j-\frac{\e^{2}}{3}\sum^{N}_{i,j=1}g(R_{p}(E_i,X)E_j, E_i)\partial_j\\
&-\frac{\e^{3}}{6}\sum^{N}_{i,j=1}(g(\n_{E_i}R_p(E_i,X)E_j,X)+g(\n_XR_p(E_i,X)E_j,E_i))\partial_{j}+P_{\e}^2(v)\sum^{N}_{j=1}\partial_{j}.
\end{align*}
Therefore
\begin{align*}
(2)+(3)&=-2\sum^{N}_{j=1}\rho_{j}\partial_j-\Delta\rho
\sum^{N}_{j=1}x^j\partial_j+\frac{2\e^{2}}{3}\sum^{N}_{i,j=1}g(R_p(E_i,E_j)E_i,X)\partial_j\\
&+\frac{\e^{3}}{3}\sum^{N}_{i,j=1}(g(\n_XR_p(E_i,X)E_i,E_j)\partial_{j}-\frac{\e^{3}}{6}\sum^{N}_{i,j=1}g(\n_{E_i}R_p(E_i,X)E_j,X)\partial_{j}\\
&+\frac{\e^{3}}{12}\sum^{N}_{i,j=1}g(\n_{E_j}R_p(E_i,X)E_i,X)\partial_{j}+P_{\e}^2(\bar{v})\sum^{N}_{j=1}\partial_{j}.
\end{align*}
Since
\begin{align*}
(1)&=\rho^{-2}\Delta-2\sum^{N}_{i,j=1}x^i\rho_{j}\partial^{2}_{ij}-\frac{\e^{2}}{3}\sum^{N}_{i,j=1}g(R_p(E_i,X)E_j,X)\partial^{2}_{ij}\\
&-\frac{\e^{2}}{6}\sum^{N}_{i,j=1}(g(\n_{X}R_p(E_i,X)E_j,X)\partial^{2}_{ij}+P_{\e}^2(\bar{v})\sum^{N}_{i,j=1}\partial^{2}_{ij},
\end{align*}
we conclude that
\begin{align*}
\Delta_{\hat{g}}&=\rho^{-2}\Delta-2\sum^{N}_{i,j=1}x^i\rho_{j}\partial^{2}_{ij}-2\sum^{N}_{j=1}\rho_{j}\partial_j-\Delta\rho
\sum^{N}_{j=1}x^j\partial_j-\frac{\e^{2}}{3} \sum^{N}_{i,j=1}g(R_p(E_i,X) E_j,X)\partial^{2}_{ij}\\
&+\frac{2\e^{2}}{3} \sum^{N}_{i,j=1}g(R_p(E_i,E_j)
E_i,X)\partial_j+\frac{\e^{3}}{3}\sum^{N}_{i,j=1}g(\n_XR_p(E_i,X)
E_i,E_j)\partial_{j}\\
&-\frac{\e^{3}}{6}\sum^{N}_{i,j=1}g(\n_{E_i}R_p(E_i,X)
E_j,X)\partial_{j}+\frac{\e^{3}}{12}\sum^{N}_{i,j=1}g(\n_{E_j}R_p(E_i,X)E_i,X)\partial_{j}\\
&-\frac{\e^{3}}{6}\sum^{N}_{i,j=1}g(\n_XR_p(E_i,X)
E_j,X)\partial^{2}_{ij}+\sum^{N}_{i,j=1}\Delta_{\e,\bar{v}}^{ij}
\end{align*}
as desired.
\QED
\section{Construction of solutions to over-determined problem}\label{s:Construction}
As explained in the previous section, our aim is to find a point $p$
and a (small) function $v:S^{N-1}\longrightarrow (0, \infty)$ such
that the over-determined problem
\begin{equation}\label{eq:problem000}
\left\{ \begin{array}{ll} -\Delta_{\bar g} \bar{u}=1& ~~\textrm{in} ~~B^{\bar{g}}_{1+v}(p)\vspace{3mm}\\
\bar{u}=0 &~~\textrm{on}  ~~\partial B^{\bar{g}}_{1+v}(p)\vspace{3mm}\\
\bar{g}(\nabla_{\bar{g}}\bar{u},\bar{\nu})=\bar{C}_0&~~\textrm{on} ~~\partial
B^{\bar{g}}_{1+v}(p),
\end{array} \right.
\end{equation}
has a solution provided $\e$ is small.
 In $\R^N$  a solution is given by
$$\phi_0(x):=\frac{1-|x|^2}{2N}$$ which clearly satisfies
 \begin{align}\label{eq:ropp2}
  \begin{cases}
    -\Delta\phi_0=1& \quad \textrm{ in }  B_1\vspace{3mm}\\
\phi_0=0&  \quad\textrm{ on }\partial B_1\vspace{3mm}\\
\dfrac{\partial \phi_0}{\partial\nu }=-\frac{1}{N}&   \quad \textrm{
on }\partial B_1.
  \end{cases}
  \end{align}
The next result  provides a first approximate solution to \eqref{eq:problem000} by solving only the Dirichlet problem in \eqref{eq:problem000}.

\begin{Proposition}\label{Propo3.2}
  There exists $\e_0>0$ such that for
all $\e \in(0,\e_0)$, for $p\in\M$  and for all $(v_0,\bar{v})\in
\mathbb{R}\times C^{2,\alpha}(S^{N-1})$ satisfying
$$|v_0|\leq\e_0,~~||\bar{v}||_{C^{2,\alpha}(S^{n-1})}\leq\e_0 \quad \textrm{ and } \quad \int_{S^{N-1}}\bar{v}\,dvol_{S^{N-1}}=0,$$
 there exists a unique positive function $\bar{\phi}=\bar{\phi}(p,\e,v_0,\bar{v}) \in C^{2,\alpha}(\ov{B^{\bar{g}}_{1+v}(p)})$ such that
\begin{equation}\label{eq:prop1}
\left\{ \begin{array}{ll} -\Delta_{\bar{g}} \bar{\phi}=1& ~~in
~B^{\bar{g}}_{1+v}(p)\vspace{3mm}\\
\bar{\phi}=0& ~~on ~~\partial B^{\bar{g}}_{1+v}(p).
\end{array} \right.
\end{equation}
 The function $\bar{\phi}$ depends smoothly on $v_0$, $\bar{v}$,  $\e$. In addition $ \bar{\phi}=\phi_0$ when $\e=0$, $v_0=0$ and  $\bar{v}\equiv0$.
\end{Proposition}
\proof  By change of variables, \eqref{eq:prop1} is equivalent to
\begin{equation}\label{eq:propp2} \left\{
\begin{array}{ll} -\Delta_{\hat{g}} \hat{\phi}=1& ~~\textrm{in}
~B_{1}\vspace{3mm}\\

\hat{\phi}=0& ~~\textrm{on}~~\partial B_{1},
\end{array} \right.
\end{equation}
where $\hat{g}$ is the induced metric defined in Lemma \ref{Lem2}.\\

Observe that, when $\e=0$, $v_0=0$ and $\bar{v}\equiv0$,
$\hat{g}$ is the Euclidean metric $g_0$ and the solution of
\eqref{eq:propp2} is given by $\hat{\phi}=\phi_0$. In fact the
solution of \eqref{eq:propp2} is the pull-back of the solution of
\eqref{eq:prop1} via the parametrization $Y_{p,v}$. We mean by this, $\hat{\phi}=Y_{p,v}^{*}\bar{\phi}$.\\
Define the Banach spaces
$$
C_{Dir}^{2,\alpha}(\ov{B_1}) :=\{ u\in C^{2,\a}(\ov{B_1})\quad:\quad
u=0\quad \textrm{ on } \de B_1\}
$$
and
$$
C^{2,\alpha}_m(S^{N-1}):=\left\{ v\in C^{2,\alpha}_m(S^{N-1}),\:\, \int_{S^{N-1}} v\, \textrm{dvol}_{S^{N-1}}=0 \right\}.
$$
Now  consider the map
\begin{align*}
 \mathcal{N}:[0,\infty)\times \mathbb{R}\times C^{2,\alpha}_m(S^{N-1})\times C_{Dir}^{2,\alpha}(B_1) &\longrightarrow C^{0,\alpha}(\ov{B_1})\\
        (\e,v_0,\bar{v},\phi) &\longmapsto\Delta_{\hat{g}}\phi+1,
 \end{align*}
It is clear that $$ \mathcal{N}(0,0,0,\phi_0)=0$$ and  $
\mathcal{N}$ is  a smooth map in a neighborhood of $(0,0,0,\phi_0)$
in $[0,\infty)\times \mathbb{R}\times C_m^{2,\alpha}(S^{N-1})\times
C_{Dir}^{2,\alpha}(\ov{B_1})$. Now since $\partial_{\phi}
\mathcal{N}(0,0,0,\phi_0)= \D: C_{Dir}^{2,\alpha}(\ov{B_1})\to
C^{0,\alpha}(\ov{B_1})$ is invertible, the implicit function theorem
gives the desired result. \QED
 Our next task is to prove that it
is possible to find $(p,\e,v_0,\bar{v})$ such that
\be\label{eq:object}
\bar{g}(\nabla_{\bar{g}}\bar{\phi},\bar{\nu})=-\frac{1}{N}\quad
\textrm{ on } ~~\partial B^{\bar{g}}_{1+v}(p). \ee
 We  compute  the Taylor of $\bar{g}(\nabla_{\bar{g}}\bar{\phi},\bar{\nu})$. To this end, we need an accurate  approximation  $ \bar{\phi}$ which is given by Proposition
\ref{Propo3.2}.
We define $\hat{\phi}=\hat{\phi}(\e,p,v): B_1\to \R $ by
\begin{equation}\label{eq:psii}
\hat{\phi}(x):= \bar{\phi}(Y_{p,v}(x))=\phi_0(\rho x)+\Psi_{\e,v}(x)~~\quad \forall x\in B_1,
\end{equation}
where we recall that $\rho=1+(v_0+  \chi v)$.
%
By \eqref{eq:propp2}, the function $\Psi_{\e,v}$ satisfies
\begin{equation}\label{eq:propp6}
\left\{
\begin{array}{ll} -\Delta_{\hat{g}} \Psi_{\e,v}=1+\Delta_{\hat{g}}\phi_0(\rho
x)~~\quad \textrm{ in }
~B_{1} \vspace{3mm}\\

\Psi_{\e,v}=-\phi_0(\rho
x) ~~\quad \textrm{ on } ~~\partial B_{1}.\\
\end{array} \right.
\end{equation}
The expansion of  $\phi_0(\rho x)$ is given by
\begin{equation}\label{ex}
\phi_0(\rho x)=\phi_{0}(x)-\frac{1}{N}|x|^{2}(\rho-1)+Q^{0}_{p}(v)
\end{equation}
and we have
\begin{Lemma}\label{lem:Psi}
The function $\Psi_{\e,v}$ defined in \eqref{eq:psii}
satisfies
\begin{align}\label{eq:sysPsi}
\begin{cases}
\displaystyle -\Delta\Psi_{\e,v}=\frac{\e^2}{3N}Ric_p(X,X) -\frac{\e^3}{4N} \sum^{N}_{i=1}g(\n_XR_p(E_i,X) E_i,X)\\
\qquad\qquad \qquad\displaystyle +\frac{\e^3}{6N}\sum^{N}_{i=1}g(\n_{E_i}R_p(E_i,X) X, X)+P_{\e}^2(v)& \quad \textrm{ in }  B_1 \vspace{3mm}\\
\displaystyle\Psi_{\e,v}=\frac{1}{N}v+Q^{0}_{p}(v)&\quad
\textrm{ on } \de B_1.
\end{cases}
\end{align}
\end{Lemma}
\proof
By straightforward computations using \eqref{ex} and  the expansion
of  $\Delta_{\hat{g}}$ in Lemma \ref{Lem2}, we get for all $x\in
B_1$,
\begin{align*}
\Delta_{\hat{g}}\phi_0(
x)&=-\rho^{-2}+\frac{4}{N}\langle\nabla\rho,x\rangle+\frac{1}{N}|x|^{2}\Delta\rho+\frac{\e^2}{3N}Ric_p(X,X)\\
&-\frac{\e^3}{4N} \sum^{N}_{i=1}g(\n_XR_p(E_i,X) E_i,X)
+\frac{\e^3}{6N}\sum^{N}_{i=1}g(\n_{E_i}R_p(E_i,X) X,
X)+P_{\e}^2(v).
\end{align*} Similarly, using Lemma \ref{Lem2}, we
have
$$
-\Delta_{\hat{g}}(\frac{1}{N}|x|^{2}(\rho-1))=-\Delta(\frac{1}{N}|x|^{2}(\rho-1))+P_{\e}^2(v)
=-\frac{1}{N}|x|^{2}\Delta\rho-\frac{4}{N}\langle\nabla\rho,x\rangle-2v+P_{\e}^2(v).
$$
From the two previous inequalities and \eqref{eq:propp6}, we deduce the first equality of
\eqref{eq:sysPsi}.   Finally  using \eqref{eq:psii}, \eqref{ex} and
the fact that $\hat{\phi}$ and $\phi_0$ are equal to $0$ on
$\partial B_1$ , we obtain
$$
\Psi_{\e,v}=\frac{1}{N}v+Q^{0}_{p}(v)~\quad \textrm{ on }~~\partial B_1.
$$
\QED
\begin{Lemma}\label{Lem3.3}
At a point $x\in \de B_1$,
we have the   expansion
\begin{align*}
 \hat{g}(\nabla_{\hat{g}}\hat{\phi},\hat{\nu})_{\mid_{\partial
B_1}}=-\frac{1}{N} +\frac{1}{N}\left[(\partial_{\nu}\psi_{v})_{|_{\partial
B_1}}- {v}\right]  +  (\partial_{\nu}\psi_\e)_{|_{\partial
B_1}}+(\partial_{\nu}\G_{\e,v})_{|_{\partial
B_1}}   +P_{\e}^1(v),
\end{align*}
where
the
functions $\psi_{{v}}$,  $\psi_{\e}$ and $\G_{\e,v}$ are respectively
(unique) solution to
\begin{equation}\label{eq:pemvw}
\displaystyle\left\{ \begin{array}{ll} \Delta \psi_{  v}=0& \quad~~\textrm{ in }~~B_1\vspace{3mm} \\
\displaystyle\psi_{  v}=  {v}&\quad ~~\textrm{ on }
~~\partial B_1,
\end{array} \right.
\end{equation}
\begin{align}\label{eqo}
\begin{cases}
\displaystyle -\Delta\psi_\e=\frac{\e^2}{3N}Ric_p(X,X) -\frac{\e^3}{4N} g(\n_XR_p(E_i,X)E_i,X)\\
\qquad\qquad \qquad\displaystyle +\frac{\e^3}{6N}g(\n_XR_p(E_i,X)X,E_i)+O_p(\e^4)& \quad \textrm{ in }  B_1 \vspace{3mm}\\
\displaystyle\psi_\e=0&\quad \textrm{ on } \de B_1.
\end{cases}
\end{align}
and
\be\label{eq:Gev}
\displaystyle\left\{ \begin{array}{ll} \Delta \G_{\e,v}= P_{\e}^2(v)& \quad~~\textrm{ in }~~B_1\vspace{3mm} \\
\displaystyle \G_{\e,v}=  Q^0(v)&\quad ~~\textrm{ on }
~~\partial B_1,
\end{array} \right.
\ee
\end{Lemma}
\proof
Since $\hat{\phi}=0$ on $B_1$, the  unit outward vector $\hat{\nu}$  about $\partial B_1$ is
given by
$$\hat{\nu}=-\frac{\nabla_{\hat{g}}\hat{\phi}}{|\nabla_{\hat{g}}\hat{\phi}|_{\hat{g}}}$$
and thus $$\hat{g}(\nabla_{\hat{g}}\hat{\phi},\hat{\nu})=-|\nabla_{\hat{g}}\hat{\phi}|_{\hat{g}}.$$
From the expansion of $\hat{g}$ in Lemma \ref{Lem2}, we have
\begin{align*}
|\nabla_{\hat{g}}\hat{\phi}|^{2}_{\hat{g}}&=\hat{g}(\nabla_{\hat{g}}\hat{\phi},\nabla_{\hat{g}}\bar{\phi})=\sum^{N}_{il=1}\hat{g}^{il}(x)\frac{\partial
\hat{\phi}}{\partial x^{i}}\frac{\partial \hat{\phi}}{\partial
x^{l}}+P_{\e}^1(v)=\rho^{-2} \sum^{N}_{i=1}\biggl(\frac{\partial
\hat{\phi}}{\partial
x^{i}}\biggl)^{2}+P_{\e}^1(v)\\
&=\rho^{-2}|\nabla\hat{\phi}|^{2}+P_{\e}^1(v).
\end{align*}
We also have $$\hat{\phi}=\phi_0(\rho x)+\Psi_{\e,v}\quad ~~
\textrm{ and }~~\quad \partial_j\phi_0(\rho
x)=-\frac{1}{N}(|x|^{2}\rho_j+x^j\rho^{2}),~j=1,...,N.$$ This
implies
$$\hat{g}(\nabla_{\hat{g}}\hat{\phi},\hat{\nu})_{\mid_{\partial
B_1}}=-\rho^{-1}|\n \hat{\phi}| +P_{\e}^1(v)=-\frac{1}{N}-\frac{1}{N}v+(\partial_{\nu}\Psi_{\e,v})_{|_{\partial
B_1}}+P_{\e}^1(v).$$
 Recalling Lemma \ref{lem:Psi}, we can decompose
$\Psi_{\e,v}$ as
\begin{equation} \label{aaa}
\Psi_{\e,v}= \frac{1}{N}\psi_{
v}+\psi_{\e}+    \G_{\e,v} ,
\end{equation}
 where the
functions $\psi_{\bar{v}}$, $\psi_{\e}$ and $\G_{\e,v}$ are respectively
(unique) solution of
$$
\displaystyle\left\{ \begin{array}{ll} \Delta \psi_{  v}=0& \quad~~\textrm{ in }~~B_1\vspace{3mm} \\
\displaystyle\psi_{  v}=  {v}&\quad ~~\textrm{ on }
~~\partial B_1,
\end{array} \right.
$$
\begin{align}\label{eqo}
\begin{cases}
\displaystyle -\Delta\psi_\e=\frac{\e^2}{3N}Ric_p(X,X) -\frac{\e^3}{4N} g(\n_XR_p(E_i,X)E_i,X)\\
\qquad\qquad \qquad\displaystyle +\frac{\e^3}{6N}g(\n_XR_p(E_i,X)X,E_i)+O_p(\e^4)& \quad \textrm{ in }  B_1 \vspace{3mm}\\
\displaystyle\psi_\e=0&\quad \textrm{ on } \de B_1
\end{cases}
\end{align}
and
$$
\displaystyle\left\{ \begin{array}{ll} \Delta \G_{\e,v}= P_{\e}^2(v)& \quad~~\textrm{ in }~~B_1\vspace{3mm} \\
\displaystyle \G_{\e,v}=  Q^0(v)&\quad ~~\textrm{ on }
~~\partial B_1.
\end{array} \right.
$$
\QED
 We
define
\begin{equation}\label{eationnn}
G(p,\e,v_0,\bar{v}):= \frac{1}{N}\left[(\partial_{\nu}\psi_{v})_{|_{\partial
B_1}}- {v}\right]  +  (\partial_{\nu}\psi_\e)_{|_{\partial
B_1}}+(\partial_{\nu}\G_{\e,v})_{|_{\partial
B_1}}   +P_{\e}^1(v),
\end{equation}
so that,\begin{equation}\label{eequationn}
\bar{g}(\nabla_{\bar{g}}\bar{\phi},\bar{\nu})|_{\de B^{\bar g}_{1+v}}=\hat{g}(\nabla_{\hat{g}}\hat{\phi},\hat{\nu})_{\mid_{\partial
B_1}}=-\frac{1}{N}+G(p,\e,v_0,\bar{v})
\end{equation}
and thus our objective \eqref{eq:object} then becomes to  find $(p,\e,v_0,\bar{v})$ such that
$G(p,\e,v_0,\bar{v})=0$. \\

 To solve this, we will use variational
 perturbative methods keeping in mind that the sets we are looking for    are stationary sets for some energy functional.
 The main strategy consists first  in  using
a local inversion argument to
  reduce the problem to finite dimensional critical point problem.
  This is due to the fact that the problem under study is invariant by
  translations on $\R^N$ and so the energy has a "kernel" at least of dimension $N$.
  This phenomenon brings some difficulties to invert the map   $ \bar v\mapsto \frac{\de}{\de  \bar v}G(p,0, 0,0)$
  as it might have    zero eigenvalues.  However, as we shall see,  $Ker \frac{\de}{\de  \bar v}G(p,0, 0,0) =\{x^i,\quad i=1,\dots,N\}$. Therefore we will solve \eqref{eq:object} modulo this set by local inversion theorems. This is the aim of  the next section.
\subsection{Local inversion argument} \label{ss:fixedpoint}
Let us consider the Dirichlet-to-Neumann operator
$$
v\mapsto (\partial_{\nu}\psi_{v})_{|_{\partial
B_1}},
$$
where
$$
\left\{ \begin{array}{ll} \Delta \psi_{ v}=0 ~~\quad \textrm{ in } ~~B_1 \vspace{3mm} \\
\psi_{  v}= {v} ~~\quad \textrm{ on } ~~\partial B_1.\\
\end{array} \right.
$$
It is well known, see for instance \cite{RaSa},  that this map has a discrete spectrum in $L^2(S^{N-1})$ given by
$$
\l_k=k,\quad k\in\N
$$
which corresponds to the Steklov eigenvalue problem. The eigenvectors corresponding to the eigenvalue $\l_k$ are given by the spherical harmonics $Y_k$ which satisfy $-\D_{S^{N-1}}Y_k=k(k+N-2) Y_k $ on $S^{N-1}$. Therefore the eigenspaces corresponding to $\l_0=0$ and $\l_1=1$ are
\be\label{eq:defE1}
\textrm{$\Lambda_0:=\textrm{span}\,\{1\}\quad$ and  $\quad \Lambda_1:=\textrm{span}\,\{x^1,\cdots,x^N\}$ }
\ee
respectively. We denote by $\Pi_0$ and $\Pi_1$ the $L^2$ projections onto  these spaces respectively and we define
$$\Pi:=\textrm{Id}-\Pi_1-\Pi_0   \qquad\textrm{and}\qquad  \Pi_1^{\perp}:=\Pi_0+\Pi. $$
 Combining these with elliptic regularity theory,  we have the following
\begin{Proposition}\label{Propo3.4}
We define the operator  $\mbL(v):=(\partial_{\nu}\psi_{v})_{|_{\partial
B_1}}- {v} $. Then
  $$\mbL:C ^{2,\alpha}(S^{N-1})\longrightarrow
C^{1,\alpha}(S^{N-1})$$
  is a self adjoint, first order elliptic
operator. In addition $$ Ker \mbL=\{x^i,\quad i=1,\dots,N\}. $$
Moreover there exists $c>0$ such that
\begin{equation}\label{L0}
||w||_{C^{2,\alpha}(S^{N-1})}\leq c
||\mbL(w)||_{C^{1,\alpha}(S^{N-1})}
\end{equation}
for every $w\in \Pi_1^\perp C^{2,\alpha}(S^{N-1})$ .
\end{Proposition}
%
We are now able to prove that, for $\e$ small enough, it is
possible to solve equation
$$
G(p,\e,v_0, \bar{v})=0
$$
modulo the kernel of $\mbL=\frac{\de}{\de    v}G(p,0, 0,0)$. Indeed we have
\begin{Proposition}\label{Pro3.5}
There exists $\e_0>0$ such that, for all
$\e\in(0,\e_0)$ and  for all $p\in \M$
 there exists a
unique $v^{\e,p}\in C^{2,\alpha}(S^{N-1}) $ with
$$
\|v^{\e,p}\|_{C^{2,\a}(S^{N-1})}<\e_0
$$
  such that
$\hat{\phi}^{\e,p}=\hat{\phi}(\e,p,v_0^{\e,p},\bar{v}^{\e,p})$ satisfies
\begin{equation}\label{eq:problem-aabb}
\left\{ \begin{array}{ll} -\Delta_{\hat g} \hat{\phi}^{\e,p}=1& ~~\textrm{in} ~~B_1 \vspace{3mm} \\
\hat{\phi}^{\e,p}=0 &~~\textrm{on}  ~~\partial B_1 \vspace{3mm}\\
\hat{g}(\nabla_{\hat{g}}\hat{\phi}^{\e,p},\hat{\nu})=-\frac{1}{N}-
\la a^{\e,p},x\ra&~~\textrm{on} ~~\partial B_1,
\end{array} \right.
\end{equation}
where $v_0^{\e,p}=\Pi_0 v^{\e,p}$,   $\la a^{\e,p},x\ra=\Pi_1 v^{\e,p}$ and   $\bar{v}^{\e,p}=\Pi v^{\e,p}$.\\

In addition  the mapping   $(\e, p,x)\mapsto {v}^{\e,p}(x)$  is smooth and satisfies
\be\label{eq:vobvae2}
  \|\n^k_g  {v}^{\e,p}\|_{ C^{2,\a}(T\M \times S^{N-1})}\leq c_k\e^{2},
\ee
for all $k\in\N $.
\end{Proposition}
\proof
We consider the map
 $$
  \calG: \M\times[0,\infty[  \times C^{2,\alpha}(S^{N-1})  \to     C^{1,\alpha}(S^{n-1})
$$
given by
$$
\calG(p,\e,v )=    G(p,\e,\Pi_0 v  ,\Pi{v}  )+ \Pi_1 v.
$$
Direct computations show that
$$
\frac{\de\calG }{\de v}(p,0,0 )[w]=
\frac{1}{N}\,  \mbL\circ \Pi_1^\perp( w)+ \Pi_1 w.
$$
We define
$$
\calL:=\frac{1}{N}\,   \mbL\circ \Pi_1^\perp+\Pi_1
$$
Thanks to Proposition \ref{Propo3.4}, the operator
$$
\calL :
C^{2,\alpha}(S^{N-1})    \to     C^{1,\alpha}(S^{N-1})
$$
 is an isomorphism and for all $w\in C^{2,\a}(S^{N-1})$
\be\label{eq:mbLpPi1}
\|w\|_{C^{2,\a}(S^{N-1})}\leq c \|\calL( w)\|_{C^{1,\a}(S^{N-1})}
\ee
Hence
 the
implicit function  theorem ensures that there exists $\e_0>0$ such that for all
$\e\in (0,\e_0)$  and for all $p\in\M$, the existence of a unique  $v^{\e,p}   \in   C^{2,\alpha}(S^{N-1})
$
  with
$$
\|v^{\e,p}\|_{C^{2,\alpha}(S^{N-1})}<\e_0
$$
 such that
$$
\calG(p,\e,v^{\e,p} )=   G(p,\e,\Pi_0 v^{\e,p} ,\Pi {v}^{\e,p}  )+\Pi_1 v^{\e,p}=0.
$$
Recalling \eqref{eationnn}, this is clearly   equivalent to
\be \label{eq:claG01}
 \calL(v^{\e,p}) +  (\partial_{\nu}\psi_\e)_{|_{\partial
B_1}}+(\partial_{\nu}\G_{\e,v^{\e,p}})_{|_{\partial
B_1}}   +P_{\e}^1(v^{\e,p}) =0.
\ee
By elliptic regularity theory
$$
\|\G_{\e,v^{\e,p}}\|_{C^{2,\a}(S^{N-1})}\leq C\e^4+ C\e^2\|v^{\e,p}\|_{C^{2,\a}(S^{N-1})}+C\|v^{\e,p}\|_{C^{2,\a}(S^{N-1})}^2.
$$
Decreasing $\e_0$ if necessary, we deduce from \eqref{eq:mbLpPi1} and \eqref{eq:claG01} that
$$  \|{v}^{\e,p}\|_{C^{2,\alpha}(S^{N-1})}\leq c\e^{2}.$$
 The smooth dependence on $\e,p$ is a consequence of the implicit function theorem. Also \eqref{eq:vobvae2} is a consequence of the fact that
$v^{\e,p}$  solves the  differential equation \eqref{eq:claG01} (which can be differentiated  $k$ times with respect to $p$) and the smooth dependence of the metric $\hat{g}$ with respect to $p$ and $\e$.
\QED
\subsection{Geometric variational argument} \label{ss:Geometric-var}
Let $\Omega\subset \M$ be a smooth   bounded domain of $\M$. It is
very well known that the minimization problem
\begin{equation}\label{eq:propp999}
J(\Omega):=\inf\biggl\{\int_{\Omega}|\nabla_{g}u|^{2}\,\textrm{
dvol}_g:~~\int_{\Omega}u\,\textrm{ dvol}_g=1,~u\in
H^{1}_{0}(\Omega)\biggl\}.
\end{equation}
has a unique solution  $u_{\Omega}\in H^{1}_{0}(\Omega)$ where
$J(\Omega)$ is achieved and  we have

\begin{equation}\label{eq:propp9}
\left\{ \begin{array}{ll}
-\Delta u_{\Omega}=J(\Omega)&\quad \textrm{ in } \Omega \vspace{3mm}\\
u_{\Omega}=0&\quad \textrm{ on } \partial \Omega.
\end{array} \right.
\end{equation}
 We can now consider the functional $\Omega\longmapsto
J(\Omega)$, for every bounded and smooth domain $\Omega\subset \M$.

\begin{Definition}
We say that $\{\Omega_{s}\}_{s\in[0,s_0)}$ is a deformation of
$\Omega_{0}$, if there exists a vector field $\Xi$ such that
$\Omega_s=\xi(s,\Omega_{0})$, where $\xi(s,.)$  is the flow
associated to $\Xi$, namely
$$\xi(0,.)=p \qquad\textrm{ and } \qquad   \frac{d\xi}{ds}(s,p)=\Xi(\xi(s,p)).$$
The deformation is volume preserving if $|\O_s|_g=|\O_0|_g$ for all $t\in[0,s_0) $.
\end{Definition}

Let $\{\Omega_{s}\}_{s\in [0,s_0)}$ be a deformation of a domain
$\Omega_{0}$ generated by the vector field $\Xi$. We denote by
$J_s=J(\O_s)$ Dirichlet's energy define in \eqref{eq:propp999},
$u_{s}$ the corresponding minimizer on $\Omega_{s}$ and  $\nu_s$ the
outward unit vector field about $\partial\Omega_s$. We have the
following lemma.

\begin{Lemma}\label{Lem2.3}
The derivative of $s\longmapsto J_s$ at $s=0$ is given by
$$
\frac{dJ_s}{ds}_{\mid_{s=0}}=-\int_{\partial \Omega_{0}}[g(\nabla_g
u_{0},\nu_{0})]^{2}g(\Xi,\nu_0) d\sigma_g, $$ where $d\sigma_g$ is
the volume element on $\partial \Omega_{0}$ for the metric induced
by $g$ and $\nu_0$ the normal vector field about $\partial
\Omega_{0}$. The domain $\Omega_0 $ is said a stationary  set for $J$ if
$$\frac{dJ_s}{ds}_{\mid_{s=0}}=0.$$
\end{Lemma}

\proof We differentiate
 \begin{equation}\label{eq:prokl}
 -\Delta_{g}u_s=J_s \quad \textrm{ in } \quad \Omega_s
\end{equation} with respect to $s$ and evaluate the result at $s=0$  to obtain
\begin{equation}\label{eq:proo}
-\Delta_{g}\partial_s u_0=J'_0  \quad \textrm{ in } \quad \Omega_0.
\end{equation}
 We also  know that
\begin{equation}\label{eq;u_s}
\int_{\Omega_s}u_{s}\textrm{ dvol}_g=1, \quad \textrm{ for  all }
~s\in [0,s_0).
\end{equation}
Differentiating  \eqref{eq;u_s} with respect to $s$ and evaluating
at $s=0$ yields
\begin{equation}\label{eq;u_tt}
\int_{\Omega_s}\partial_su_0\textrm{ dvol}_g=0.
\end{equation}
We multiply \eqref{eq:proo} by $u_{0}$ and \eqref{eq:prokl},
evaluated at $s=0$, by $\partial_su_0$, subtract the two results and
integrate over $\Omega_0$ to get
 \begin{align*}
J'_0\int_{\Omega_0}u_{0}\textrm{
dvol}_g-J_0\int_{\Omega_0}\partial_su_{0}\textrm{
dvol}_g&=\int_{\Omega_0}\biggl(\partial_su_0\Delta_{g}
u_{0}-u_{0}\Delta_{g}\partial_su_0\biggl)\textrm{
dvol}_g\\
&=\int_{\partial\Omega_0}\biggl(\partial_su_0g(\nabla_g u_0,\nu_{0})-u_{0}\frac{\partial(\partial_su_0)}{\partial\nu_0}\biggl)d\sigma_g\\
&=\int_{\partial\Omega_0}\partial_su_0g(\nabla_g
u_0,\nu_{0})d\sigma_g,
\end{align*}
where we have used the fact that $u_0=0~~on~\partial\Omega_0$ to
obtain the last equality. We conclude with \eqref{eq;u_s} and
\eqref{eq;u_tt} that
\begin{equation}\label{j'}
J'_0= \int_{\partial\Omega_0}\partial_su_0g(\nabla_g
u_0,\nu_{0})d\sigma_g
 \end{equation}
Now, let  $\xi$ be the flow generated by $\Xi$, by definition
\begin{equation}\label{eq;u_ttt}
u_s(\xi(s,p))=0~~ \quad\textrm{ for } ~~p\in
\partial \Omega_s.
~\end{equation} We differentiate \eqref{eq;u_ttt} with respect to
$t$ and evaluating at $s=0$ and get

$$\partial_su_0=-g(\nabla_g u_0,\Xi).$$

But, $u_0=0$ on $\partial \Omega_0$, and hence only the normal
component of $\Xi$ plays a role in this formula. Therefore, we have
$$
\partial_su_0=-g(\nabla_g u_0,\nu_{0})\, g(\Xi,\nu_0)~on~\partial
\Omega_0
$$
 and replacing this in \eqref{j'}, we finally get that
$$
J'_0=-\int_{\partial \Omega_0}[g(\nabla_g u_0,\nu_0)]^{2}\,g(\Xi,\nu_{0})\, d\sigma_g.
$$
 \QED
 The following proposition gives a necessary  and sufficient
condition for a domain $\O$ being a stationary set of $J$.
\begin{Proposition}
A domain $\O$ is a stationary set for $J$ under volume preserving deformations if and only if their exits a function $u_{\O_0}$ such that
\begin{align}\label{eq:problembJ}
  \begin{cases}
\displaystyle    -\Delta_{g} u_{\O_0}=J(\O_0)& \quad \textrm{ in }  \Omega_0 \vspace{3mm}\\
\displaystyle  u_{\O_0}=0&  \quad\textrm{ on }\partial\Omega_0 \vspace{3mm}\\
\displaystyle  g(\n_g u_{\O_0}, \nu_0)=\l &   \quad \textrm{ on }\partial\Omega_0,
  \end{cases}
  \end{align}
  for some $\l\in\R$.
\end{Proposition}
The  proof of this proposition is similar to the one  of  [Proposition
2.2 in \cite{FPP}] so we skip it.\\

We also remark that instead of considering volume preserving deformation, a smooth bounded stationary set $\O_0$ for the total energy
$$
\O\mapsto J(\O)+\l^2 |\O|_g
$$
implies the existence of $u_{\O_0}$ such that  \eqref{eq:problembJ} holds. This can be seen from
Lemma \ref{Lem2.3} and the variation of volume which is given by
$$
\frac{d}{ds}|\O_s|_g |_{s=0}=\int_{\de \O_0}g(\Xi,\nu_{0})\, d\sigma_g.
$$
See for instance [\cite{He}, Theorem 1.11].
\subsubsection{The reduced functional}
Let us recall what we have obtained so far. Thanks to Proposition \ref{Pro3.5} and the usual $Y_{p,v^{\e,p}}$-change of variable,
we have: for all $p\in\M$ and for all $\e>0$ small we have  $\bar{\phi}^{\e,p} =\hat{\phi}^{\e,p}\circ Y^{-1}_{p,v^{\e,p}}$
satisfies
\begin{equation}\label{eq:problem-a}
\left\{ \begin{array}{ll} -\Delta_{\bar g} \bar{\phi}^{\e,p}=1& ~~\textrm{in} ~~B^{\bar{g}}_{1+v^{\e,p}}(p)\vspace{3mm}\\
\bar{\phi}^{\e,p}=0 &~~\textrm{on}  ~~\partial B^{\bar{g}}_{1+v^{\e,p}}(p)\vspace{3mm}\\
\bar{g}(\nabla_{\bar{g}}\bar{\phi}^{\e,p},\bar{\nu})=-\frac{1}{N}-\bar{g}(A^{\e,p},{\calV}^{\e,p})&~~\textrm{on}
~~\partial B^{\bar{g}}_{1+v^{\e,p}}(p),
\end{array} \right.
\end{equation}
where for all $x\in S^{N-1}$, we define $
A^{\e,p}(Y_{p,v^{\e,p}}(x)):=dY_{p,v^{\e,p}}(x)[ a^{\e,p}] $ and
similarly    ${\calV}^{\e,p}:=dY_{p,v^{\e,p}}(x)[ x]  $.
It follows that the inverse of the torsion rigidity for $ \bar{\phi}^{\e,p}$ is given by
$$
J( B^{\bar{g}}_{1+v^{\e,p}}(p))=\frac{1}{
\displaystyle\int_{B^{\bar{g}}_{1+v^{\e,p}}(p)}\bar{\phi}^{\e,p}\textrm{dvol}_{\bar
g}}.
$$
The domains $\O_\e$ we are looking for is a critical point of the  the total energy functional :
$$
\O \mapsto J(
\O)+\frac{1}{N^2}|\O
|_{\bar{g}}.
$$
This allows to define for $p\in \M$, the reduced  functional
\begin{equation} \label{eq:redF}
\displaystyle \Phi_\e(p):=J(
B^{\bar{g}}_{1+v^{\e,p}}(p))+\frac{1}{N^2}|B^{\bar{g}}_{1+v^{\e,p}}(p)
|_{\bar{g}}.
\end{equation}

\begin{Proposition} \label{prop:soluptoscale}
Let $\bar{\phi}^{\e,p}$ satisfies  \eqref{eq:problem-a}.  If $p$ is a critical point of $\Phi_\e$ then $A^{\e,p}=0$, provided $\e$ is small. In particular
$$
\left\{ \begin{array}{ll} -\Delta_{\bar g} \bar{\phi}^{\e,p}=1& ~~\textrm{in} ~~B^{\bar{g}}_{1+v^{\e,p}}(p)\vspace{3mm}\\
\bar{\phi}^{\e,p}=0 &~~\textrm{on}  ~~\partial B^{\bar{g}}_{1+v^{\e,p}}(p)\vspace{3mm}\\
\bar{g}(\nabla_{\bar{g}}\bar{\phi}^{\e,p},\bar{\nu})=-\frac{1}{N}
&~~\textrm{on} ~~\partial B^{\bar{g}}_{1+v^{\e,p}}(p).
\end{array} \right.
$$
\end{Proposition}

\proof Given  $\Xi\in T_{p}\M$,  we consider the geodesic curve
$p_s=\exp_p(s\Xi)$. Let $E^s_i$ be the parallel transport of $E_i$ to $p_s$ along the curve $[0,1]\ni t\mapsto \exp_{p}^{\bar g }(ts E_i)$. Provided $s$ is fixed and small, we can consider
the perturbed ball $B^{\bar{g}}_{1+v^{\e,p_s}}(p_s)$  so that
\eqref{eq:problem-a} holds.  Recall that
$$
B^{\bar{g}}_{1+v^{\e,p_s}}(p_s)=Y_{v^{\e,p_s},p_s}(B_1).
$$
Define the vector field
$$
W_s( Y_{p,v^{\e,p}}(x))=(1+
v^{\e,p_s}(x))\sum_{i=1}^Nx^i E^s_i \quad \forall x\in B_1.
$$
 We   now define the deformation of $
B^{\bar{g}}_{1+v^{\e,p}}(p) $ by
$$
\xi(s, q)=  \exp_{p_s}(W_s(q))\quad \forall q \in
Y_{p,v^{\e,p}}(B_1).
$$
Next we observe that
$$
\frac{d \xi}{ds }(0,q)=J_q(1),
$$
  where $J_q(t)=\de_s\exp_{p_s}(tW_s(q))|_{s=0}$ is the  Jacobi field along the geodesic
$\g_q(t)= \exp_{p}(tW_0(q))$ with
$$
\textrm{ $J_q(0)=\Xi\qquad $ and $\qquad J_q'(0):=\frac{D J_q}{dt}(0)= \frac{D W_s(q)}{ds}|_{s=0}.$}
$$
Note that for   $q= Y_{p,v^{\e,p}}(x)$, we have
$$
J_q'(0)= \frac{D W_s(q)}{ds}|_{s=0} = d_p v^{\e,p}(x)[\Xi]X+(1+v^{\e,p}(x))x^i \frac{D E^s_i}{ds}|_{s=0}=d_p v^{\e,p}(x)[\Xi]X
$$
and thus by \eqref{eq:vobvae2}, we get
\be\label{eq:Jqp0e2}
|J_q'(0) |\leq C\e^2 |\Xi|_g.
\ee
Since also  $t\mapsto J_q(t)$ satisfies an homogenous  second order linear differential
equation with uniformly bounded coefficients with respect to $\e$ and $q$, we get for all $q\in Y_{p,v^{\e,p}}(B_1)$ \be
\label{eq:estDjdtJ1} |J_q(1)|_{\bar g}\leq C ( |J_q'(0)|_g +  | J_q(0)|_g )\leq C |\Xi|_g. \ee Thanks
to  [Proposition 3.6, in \cite{doCarmo}], we have \be \bar
g(J_q(1),\g_q'(1))=\bar g\left( J_q'(0),\g_q'(0)\right)+ \bar g(
J_q(0),\g_q'(0)). \ee
It is plain that  at any  point $q=Y_{p,v^{\e,p}}(x) \in \de  B^{\bar{g}}_{1+v^{\e,p}}(p) $
$$
\bar g(J_q(1),\bar{\nu}(q)) = \bar g(J_q(1),\bar{\nu}(q)- \g_q'(1))+ \bar g(J_q(1),  \g_q'(1))
$$
which implies
\begin{align*}
\bar g(J_q(1),\bar{\nu}(q)) -\bar g(\Xi,X) &= \bar g(J_q(1),\bar{\nu}(q)- \g_q'(1))+ \bar g\left( J_q'(0),\g_q'(0)\right)+ \bar g( \Xi,\g_q'(0) -X) \\
&= \bar g(J_q(1),\bar{\nu}(q)- \g_q'(1))+ \bar g\left( J_q'(0),(1+v^{\e,p}) X\right)+v^{\e,p} \bar g( \Xi,  X ),
\end{align*}
where we have  used the fact that $ \g_q'(0)=W_0(q)=(1+v^{\e,p}(x))
X$. We also have (see for instance \cite{PX} for the expansion of
$-\bar{\nu}(q)$ and recall \eqref{eq:Upsilon})
$$
|\bar{\nu}(q)- \g_q'(1)|_{\bar g}= |\bar{\nu}(q)- (1+v^{\e,p}(x))\Upsilon (x)|_{\bar g} \leq C \e^2.
$$
By using this \eqref{eq:vobvae2}, \eqref{eq:Jqp0e2} and \eqref{eq:estDjdtJ1}, we then  deduce that, at any point $q=Y_{p,v^{\e,p}}(x)\in \de  B^{\bar{g}}_{1+v^{\e,p}}(p)  $,
\be \label{eq:ZXi}
|\bar g(J_q(1),\bar{\nu}(q))-\bar g(\Xi,X)|_g \leq C \e^2 |\Xi|_g.
\ee
We now recall that
$$
\displaystyle \Phi_\e(p_s):=J(
\xi(s,B^{\bar{g}}_{1+v^{\e,p}}(p)))+\frac{1}{N^2}|\xi(s,B^{\bar{g}}_{1+v^{\e,p}}(p))
|_{\bar{g}}.
$$
 Assume now that $p$ is a critical point of $\Phi_\e$. Then by  Lemma \ref{Lem2.3} and the last equation of \eqref{eq:problem-a}, we have
\begin{align*}
0= \frac{d}{ds}\Phi_\e(p_s)_{\mid_{s=0}}&=- \int_{\partial
B^{\bar{g}}_{1+v^{\e,p}}(p)}[\bar{g}(\nabla_{\bar{g}}\bar{\phi},\bar{\nu})]^{2}\bar{g}(J_q(1),\bar{\nu})
d\sigma_{\bar{g}}+ \frac{1}{N^2}\int_{\partial
 B^{\bar{g}}_{1+v^{\e,p}}(p)} \bar{g}(J_q(1),\bar{\nu})
d\sigma_{\bar{g}}\\
&=- \frac{2}{N}\int_{\partial B^{\bar{g}}_{1+v^{\e,p}}(p)}
\bar{g}(J_q(1),\bar{\nu}) d\sigma_{\bar{g}}- \int_{\partial
B^{\bar{g}}_{1+v^{\e,p}}(p)}
\bar{g}(A^{\e,p},{\calV}^{\e,p})^{2}\bar{g}(J_q(1),\bar{\nu})
d\sigma_{\bar{g}}.
\end{align*}
Hence
$$
\frac{2}{N}\int_{\partial  B^{\bar{g}}_{1+v^{\e,p}}(p)
}\bar{g}(A^{\e,p},{\calV}^{\e,p}) \bar{g}(J_q(1),\bar{\nu})
d\sigma_{\bar{g}}=-\int_{\partial
B^{\bar{g}}_{1+v^{\e,p}}(p)}\bar{g}(A^{\e,p},{\calV}^{\e,p})^{2}
\bar{g}(J_q(1),\bar{\nu}) d\sigma_{\bar{g}}
$$
and by \eqref{eq:estDjdtJ1}
$$
\frac{2}{N} \int_{\partial  B^{\bar{g}}_{1+v^{\e,p}}(p)
}\bar{g}(A^{\e,p},{\calV}^{\e,p}) \bar{g}(J_q(1),\bar{\nu})
d\sigma_{\bar{g}}    \leq  c  ||\Xi||_{\hat{g}} \int_{\partial
B^{\bar{g}}_{1+v^{\e,p}}(p)}\bar{g}(A^{\e,p},{\calV}^{\e,p})^{2}
  d\sigma_{\bar{g}}
$$
 for all $\Xi\in T_{p}\M$.
 By    changing   variables, using \eqref{eq:estDjdtJ1} and \eqref{eq:ZXi},   we obtain
\begin{align*}
\frac{2}{N}\int_{\partial B_1} \la a^{\e,p},x\ra
\langle\Xi,X\rangle\, d\sigma_{\hat{g}}(x)\leq
c\e^2||\Xi||_{\hat{g}}| a^{\e,p}|+ c  ||\Xi||_{\hat{g}}
\int_{\partial B_1}\la a^{\e,p},x\ra^{2}\, d\sigma_{\hat{g}}(x).
\end{align*}
From Lemma \ref{Lem2}, we get
\begin{align*}
\frac{2}{N}\int_{\partial B_1}  \la a^{\e,p},x\ra
\langle\Xi,X\rangle\, dx  \leq c\e^2||\Xi||_{\hat{g}}| a^{\e,p}|+ c
||\Xi||_{\hat{g}} \int_{\partial B_1}\la a^{\e,p},x\ra^{2}\, dx.
\end{align*}
We now choose $\Xi=\sum_{i=1}^N a^{\e,p}_i E_i\in T_p \M$ and use the
fact that $|a^{\e,p}|\leq c\e^2$ to get
$$
\frac{1}{2N}\int_{\partial B_1}\langle a^{\e,p} ,x\rangle^2 dx\leq
C\e^2 |  a^{\e,p}|^2(1+|  a^{\e,p}|),
$$
for some positive constant $C$ provided $\e$ is small.  We then
conclude that
$$
 \frac{1}{2N}|  a^{\e,p}|^2 \leq C\e^2 |  a^{\e,p}|^2(1+|  a^{\e,p}|),
$$
provided $\e$ is small enough. This shows that  $a^{\e,p}=0$.
\QED
\subsection{Expansion of volumes of the  perturbed geodesic ball}

\begin{Lemma}\label{Lem4.1}
Assume that $v_0^{\e,p}$ is given by Proposition \ref{Pro3.5}. Then for all $\e$ positive small, we have
$$v_0^{\e,p}=-\frac{S_g(p)}{3N(N+2)}\e^{2}+O_p(\e^{4}).$$
\end{Lemma}
\proof Under the hypothesis of Lemma \ref{Lem4.1}, we have with
Proposition \ref{Pro3.5} that
$$\hat{g}(\nabla_{\hat{g}}\hat{\phi}^{\e,p},\hat{\nu})=-\frac{1}{N}-
\la a^{\e,p},x\ra~~\textrm{on} ~~\partial B_1.$$ Using
\eqref{eequationn}, we get equivalently
$$G(p,\e,v_0^{\e,p},\bar{v}^{\e,p})+\langle
a^{\e,p}, x\rangle=0~~\textrm{on} ~~\partial B_1.$$ Equation
\eqref{aaa} together with the estimate in proposition \ref{Pro3.5}
yield
\begin{equation}\label{eqv0}
\frac{1}{N}\mbL(\bar{v}^{\e,p})-\frac{1}{N}v_0^{\e,p}+(\partial_{\nu}\psi_{\e})_{\mid_{\partial
B_{1}}}+\langle a^{\e,p},x\rangle+O_p(\e^{4})=0,
\end{equation}
where  $\psi_{\e}$ is solution of  \eqref{eqo}. Because the
integral of the maps $ \mbL(\bar{v}^{\e,p})$  and $\langle
a^{\e,p},x\rangle$ over $S^{N-1}$ are equal to $0$, we get
integrating \eqref{eqv0} that
\begin{align*}
v_0^{\e,p}|B_1|&=\int_{S^{N-1}}\partial_{\nu}\psi_{\e}\, \textrm{dvol}_{S^{N-1}}+O_p(\e^{4})\\
&=\int_{B_1}\Delta\psi_{\e}+O_p(\e^{4})=\frac{1}{3N}\sum^{N}_{i,j,k=1}R_{ijik}\int_{B_1}x^{k}x^j\e^{2}+O_p(\e^{4}),
\end{align*}
where we have used \eqref{eqo}, the notation
$$
R_{ijkl}=g(R_{p}(E_i,E_k)E_j,E_l)~\quad \textrm{ and } \quad ~R_{ikjl,m}=g(\nabla_{E_{m}}R_{p}(E_i,E_k)E_j,E_l)$$
and the fact that the integral of a spherical harmonic odd degree
over the unit sphere $S^{N-1}$ is equal to $0$.
 Now, using the identity
 \begin{equation}\label{eq:iddd}
 \int_{\partial
B_1}x^{k}x^{l} \, \textrm{dvol}_{S^{N-1}}=|B_1|\delta_{kl},
 \end{equation}
 we deduce that
\begin{equation}\label{eqvv0}
v_0^{\e,p}=-\frac{S_g(p)}{3N(N+2)}\e^{2}+O_p(\e^{4}),
\end{equation}
\QED

\begin{Proposition}\label{Provolume}
Assume that $v_0^{\e,p},\bar{v}^{\e,p}$ are as in Proposition \ref{Pro3.5}.
Then as $\e\to 0$,  we have
$$
|\partial
B^{g}_{\e(1+v^{\e,p})}(p)|_{g}=N|B_1|\e^{N-1}\biggl(1-\frac{N+4}{6(N+2)}S_g(p
)\e^{2}+O_{p}(\e^{4})\biggl)
$$
and
$$
|B^{g}_{\e(1+v^{\e,p})}(p)|_{g}=|B_1|\e^{N}\biggl(1-\frac{1}{2(N+2)}S_g(p
)\e^{2}+O_{p}(\e^{4})\biggl).
$$
\end{Proposition}
\proof Recall that $\bar{g}=\e^{-2}g$ and this implies
$$
|B^{g}_{\e(1+v^{\e,p})}(p)|_g=\e^{N}|B^{\bar{g}}_{1+v^{\e,p}}(p)|_{\bar{g}}=\e^{N}|B_1|_{\hat{g}}
$$~~
and~~
$$
|\partial B^{g}_{\e(1+v^{\e,p})}(p)|_g=\e^{N-1}|\partial
B_1|_{\hat{g}}.
$$

We get from the expansion in Lemma \ref{Lem2} that
\begin{align}\label{det}
\sqrt{|\hat{g}|}&=1+Nv_0^{\e,p}+N\chi\bar{v}^{\e,p}+\langle
x,\nabla\rho\rangle+\frac{1}{6}\sum^{N}_{k,l,s=1}R_{sksl}x^{k}x^{l}\e^{2} \nonumber\\
&+\frac{1}{12}\sum^{N}_{k,l,s,m=1}R_{sksl,m}x^{k}x^{l}x^m\e^{3}+O_p(\e^{4}),
\end{align}
 The expansion of $|\partial B_1|_{\hat{g}}$ then follows  integrating \eqref{det} over the unit sphere $S^{N-1}$,
where we use the value of $v^{\e,p}_0$ in Lemma \ref{Lem4.1}, the
identity \eqref{eq:iddd} and the fact that, the integral over the
unit sphere $S^{N-1}$ of a spherical harmonic of odd degree is equal
to $0$ and the function $\bar{v}^{\e,p}$ has mean value equal to
$0$. Similarly we get $|B_1|_{\hat{g}}$ by   integrating
\eqref{det} over the unit ball $B_1$. \QED

{\section{Proof of the Theorem \ref{theo0}}} In the following result, we  characterize critical points of the function
$\Phi_{\e}$ leading to the location  of the
extremal domains we have constructed in the previous sections. We
recall the reduced functional defined in \eqref{eq:redF} by
$$\Phi_\e(p):=J(
B^{\bar{g}}_{1+v^{\e,p}}(p))+\frac{1}{N^2}|B^{\bar{g}}_{1+v^{\e,p}}(p)
|_{\bar{g}}.$$
\begin{Lemma} \label{Lem4.3}
As $\e$ tends to zero, we have

\be\label{eq:phiep0} \Phi_\e(p)=\alpha_{N}+\beta_{N}\, \e^2\, S_g(p) + O_p(\e^4), \ee

where

$$\alpha_{N}=\frac{N^{3}(N+2)+|B_1|^{2}}{N^{2}|B_1|} \qquad \textrm{ and }\qquad  \beta_{N}=\frac{N^{2}(N+2)^{3}-(N+4)|B_1|^{2}}{2N^{2}(N+2)(N+4)}.$$
In addition $\b_N\neq0 $ for every $N\geq 2$.
\end{Lemma}
\proof After change of variable, we can write  $\Phi_\e$ on
the form
\be \label{eq:phie01}
 \Phi_\e(p):=\frac{1}{ \displaystyle\int_{B_1}\hat{\phi}^{\e,p}\textrm{dvol}_{\hat g}}+\frac{1}{N^2}|B_1
|_{\hat{g}} .
\ee

From the estimate of $\bar{v}^{\e,p}$ in Proposition \ref{Pro3.5},
$\hat{\phi}^{\e,p}$ is now written as
\begin{equation}\label{eq:pkk}
\hat{\phi}^{\e,p}=\phi_{0}-|x|^{2}\frac{v_0^{\e,p}}{N}-\frac{1}{N}|x|^{2}\chi\bar{v}^{\e,p}+\Psi_{\e,v^{\e,p}}+O_{p}(\e^{4}),
\end{equation}
where $\Psi_{\e,v^{\e,p}}$ is given by Lemma  \ref{lem:Psi} with
$\rho=v_0^{\e,p}+\chi\bar{v}^{\e,p} $. We integrate the function
$\hat{\phi}^{\e,p}$ over the unit ball  $B_1$ using the volume element of
$\hat{g}$.

We get, using \eqref{det} and the fact that  $\bar{v}^{\e,p}$ has  zero mean value,
\begin{align}\label{eq psi}
\int_{B_1}\hat{\phi}^{\e,p}\textrm{dvol}_{\hat
g}&=(1+Nv^{\e,p}_0)\int_{B_1}\phi_{0}-\frac{v^{\e,p}_0}{N}\int_{B_1}|x|^{2}\nonumber\\
&+\int_{B_1}\Psi_{\e,v}+\frac{1}{6}\sum^{N}_{k,l,s=1}R_{sksl}\e^{2}\int_{B_1}x^{k}x^{l}\phi_0+O_{p}(\e^{4}).
\end{align}
A straightforward computation yields
\begin{equation}\label{eq comput}
\int_{B_1}\phi_{0}(|x|)=\frac{|B_1|}{N(N+2)}, \qquad \qquad \int_{B_1}|x|^{2}=\frac{N|B_1|}{N+2}
\ee
 { and }
\be\label{eq psi1}
\sum^{N}_{k,l,s=1}R_{sksl}\e^{2}\int_{B_1}x^{k}x^{l}\phi_0\, \textrm{dvol}_{S^{N-1}}=\frac{-|B_1|S_g(p)}{N(N+2)(N+4)}\e^{2},
\end{equation}
In other to compute the integral of $\Psi_{\e,v^{\e,p}}$ over $B_1$, we use the
formula
   \begin{equation}\label{eq:pkl}
   \int_{\partial B_1}\biggl(\phi_{0}\frac{\partial \Psi_{\e,v}}{\partial
   \nu}-\Psi_{\e,v}\frac{\partial\phi_0}{\partial\nu}\biggl)\, \textrm{dvol}_{S^{N-1}} =\int_{B_1}\biggl(\phi_{0}\Delta\Psi_{\e,v}-\Psi_{\e,v}\Delta\phi_0\biggl).
   \end{equation}
Recall that  $\frac{\partial\phi_0}{\partial \nu}=-\frac{1}{N}$, $\phi_0=0$ on $\partial
   B_1$ and $-\Delta\phi_0=1$ in $B_1$. Using this, we get from \eqref{eq:pkl}
   that
\begin{align*}
\int_{B_1}\Psi_{\e,v}& =\frac{1}{N}\int_{\partial
B_1}\Psi_{\e,v}\, \textrm{dvol}_{S^{N-1}}-\int_{B_1}\phi_{0}\Delta\Psi_{\e,v}\\
&=\frac{|B_1|}{N}v^{\e,p}_0-\frac{1}{3N}\sum^{N}_{i,j,k=1}R_{ijik}\e^{2}\int_{B_1}\phi_{0}x^{k}x^j +O_{p}(\e^{4}).
\end{align*}
Therefore
\begin{equation}\label{eq intPsi}
\int_{B_1}\Psi_{\e,p}=\frac{|B_1|}{N}v^{\e,p}_0+\frac{|B_1|S_g(p)}{3N^{2}(N+2)(N+4)}\e^{2}+O_{p}(\e^{4}).
\end{equation}
Replacing \eqref{eq intPsi}, \eqref{eq comput} and \eqref{eq psi1} in \eqref{eq psi},
we obtain
\begin{equation}\label{critic}
\int_{B_1}\hat{\phi}^{\e,p}\textrm{dvol}_{\hat{g}}=\frac{|B_1|}{N(N+2)}\biggl(1+(N+2)v^{\e,p}_0-\frac{(N-2)S_g(p)}{6N(N+4)}\e^{2}+O_{p}(\e^{4})\biggl).
\end{equation}
One can now consider the value of $v^{\e,p}_0$ to get
\begin{equation}\label{criticc}
\int_{B_1}\hat{\phi}^{\e,p}\textrm{dvol}_{\hat
{g}}=\frac{|B_1|}{N(N+2)}\biggl(1-\frac{N+2}{2N(N+4)}S_g(p)\e^{2}+O_{p}(\e^{4})\biggl).
\end{equation}
That is
$$\frac{1}{ \displaystyle\int_{B^{\bar{g}}_{1+v^{\e,p}}(p)}\hat{\phi}^{\e,p}\textrm{dvol}_{\hat g}}=\biggl(1+\frac{N+2}{2N(N+4)}S_g(p)\e^{2}+O_{p}(\e^{4}))\biggl)J_1,$$
where
$$
J_1=\frac{N(N+2)}{|B_1|}.
$$
We now use the expansion of
$|B_1|_{\hat{g}}$ in Proposition \ref{Provolume} which we plug in \eqref{eq:phie01} to get \eqref{eq:phiep0}.\\
Next we prove that $\b_N\neq0$.
 Suppose on the contrary that  for some  integer $N\geq 2$   we have
$\beta_{N}=0$. Then
\begin{equation}\label{eq:beta0}
|B_1|^{2}=\frac{N^{2}(N+2)^{3}}{N+4}.
\end{equation}
We now recall  the volume of the unit ball $|B_1|$ in $\R^N$.  For  $N=2k$, an even integer, it is given by
$$
|B_1| =\frac{\pi^{k}}{2k k!}
$$
and for $N=2k+1$  we have
$$
~~~|B_1| =\frac{2^{2k+1}\pi^{k}k!}{(2k+1)(2k+1)!}.
$$
These imply that
$$
 \frac{16k^{2}(k+1)^{3}}{k+2}=\biggl(\frac{\pi^{k}}{2kk!}\bigg)^{2}
$$
 and
 $$
 \frac{(2k+1)^{2}(2k+3)^{3}}{2k+5}=\biggl(\frac{2^{2k+1}\pi^{k}k!}{(2k+1)(2k+1)!}\biggl)^{2}.
$$
The above equalities contradict the fact that $\pi$ is a  {transcendental number}, see \cite{Lind}.
\QED
We now  complete the
proof of Theorem \ref{theo0} by defining
$$
\calF(p,\e):=\frac{1}{\beta_{N}}\frac{\Phi_{\e}(p)-\alpha_{N}}{\e^{2}}
.$$ It follows from Lemma \ref{Lem4.3} that
$$
\|\calF(\cdot,\e)-S_g \|_{C^{2,\a}(\M)}\leq C\e^{2},
$$
for a positive constant $C$ independent of $\e$. If $p$ is a
critical point of $\calF(\cdot,\e) $ then by Proposition
\ref{prop:soluptoscale}, we have
$$
\left\{ \begin{array}{ll} -\Delta_{\bar g} \bar{\phi}^{\e,p}=1& ~~\textrm{in} ~~B^{\bar{g}}_{1+v^{\e,p}}(p)\vspace{3mm}\\
\bar{\phi}^{\e,p}=0 &~~\textrm{on}  ~~\partial B^{\bar{g}}_{1+v^{\e,p}}(p)\vspace{3mm}\\
\bar{g}(\nabla_{\bar{g}}\bar{\phi}^{\e,p},\bar{\nu})=-\frac{1}{N}
&~~\textrm{on} ~~\partial B^{\bar{g}}_{1+v^{\e,p}}(p_).
\end{array} \right.
$$
Now we recall that $\bar{g}=\e^{-2}g$ and so we put     $u_\e=
\e^2 \bar{\phi}^{\e,p}$ and $\O_\e= B^{ {g}}_{\e(1+v^{\e,p})}(p)$. It is also clear from the construction that
$$
\|u_\e\|_{C^{2}(\ov{ \Omega_{\e}})}\leq C.
$$
We therefore finish the proof of Theorem \ref{theo0}.
\QED
\section{Local foliation by boundaries of extremal domains}\label{s:local-fol}
Let $E_i^t$  be the parallel transport of $E_i$  along the geodesic $\exp_{p_0}(t E_i)$ for all $i=1,\dots,N$.\\
For $\t\in\R^N$, we let $q=\exp_{p_0}(\t^i E_i)$ and
consider as usual
$$
Y_{q, v}(x):=\exp_{q}^{\bar g}\left((1+v_0+\chi \bar{v})\sum_{i=1}^N x^i E_i^{\t^i}\right).
$$
Then there exists $v^{\e,q}$ such that  $B^{\bar g}_{1+v^{\e,q}}(q)=Y_{\e,v^{\e,q}}(B_1)$  satisfies \eqref{eq:problem-a}.\\

Let us now assume that  $p_0$ is a non-degenerate critical point of the scalar
curvature function $S_g$. Then by the implicit function theorem,  Proposition \ref{prop:soluptoscale} and
 Lemma \ref{Lem4.3} there exists a regular curve $\t(\e)\in \R^N$ with $|\t(\e)|\leq C \e^2$ and
such that
$$
\n_g\calF(\e, q_\e)=0,
$$
where $ q_\e= \exp_{p_0}\left(\sum_{i=1}^N\t^i(\e) E_i\right)$. Therefore  by Proposition \ref{prop:soluptoscale} and a scaling argument we have a smooth function $\phi^\e= \phi^{\e,q_\e} $ such that
$$
\left\{ \begin{array}{ll} -\Delta_{  g}  \phi^\e=1& ~~\textrm{in} ~~B^{ {g}}_{\e(1+v^{\e,q\e})}(q_\e)\vspace{3mm}\\
\phi^\e=0 &~~\textrm{on}  ~~\partial B^{ {g}}_{\e(1+v^{\e,q_\e})}(q_\e)\vspace{3mm}\\
 {g}(\nabla_{ {g}}\phi^\e, {\nu}_\e)=-\frac{\e}{N}
&~~\textrm{on} ~~\partial B^{ {g}}_{\e(1+v^{\e,q_\e})}(q_\e).
\end{array} \right.
$$

We will prove in our next result that the family of hyper-manifolds $$\left(\de B^{ {g}}_{\e(1+v^{\e,q_\e})}(q_\e),\quad\e\in (0,\e_0) \right)$$ constitutes a foliation.  This is an immediate consequence  of Proposition \ref{prop:folli} below. The main ingredients of the proof is contained in Ye \cite{Ye}. However we will
  write  a more  applicable result.
\begin{Proposition}\label{prop:folli}
Let  $p_0 \in\M$  and $\g: [0, t_0]\to \M$ be a regular curve such that $\g(0)=p_0 $ and $|\g'(0)|_g=0$. Let $v: [0,t_0]\times S^{N-1}\to \R$ be a $C^2$-function such that $v(0,\cdot) =0$.\\

Then there exists $t_1\in (0,t_0)$ and a $C^2$-function $\o : (0,t_1]\times S^{N-1}\to \R_+ ^*$ such that  for all $t\in (0,t_1]$
$$
 \left\{\exp_{\g(t)}\left( t(1+v(t,x))\sum^{N}_{i=1}x^i E_i^t \right)\,:\, x\in S^{N-1} \right\} =\left\{\exp_{p_0}\left( {\o}(t, y)\sum^{N}_{i=1}  y^i E_i\right)\,:\, y\in S^{N-1} \right\},
$$
where $E_i^t $ is     the parallel transport to $\g(t)$ of $E_i$ along the geodesic $s\mapsto \exp_{p_0}(s E_i)$.\\

In addition $$
\de_{t} {\o}(0, \cdot)=1.
 $$
 In particular setting $ S^{ {g}}_{t (1+v^t )}(\g(t))= \left\{\exp_{\g(t)}\left( t(1+v(t,x))\sum^{N}_{i=1}x^i E_i^t \right)\,:\, x\in S^{N-1} \right\} $, then  the family of  perturbed balls
 $\left( S^{ {g}}_{t (1+v^t )}(\g(t)),\quad t\in (0,t_1)\right)$ constitutes a
smooth foliation of a neighborhood of $p_0$.
\end{Proposition}
\proof To alleviate the notations, we put $v^t=v(t,\cdot)$ and $p_t=\g(t)$.\\
 \noindent \textbf{Claim:} There exists a smooth function
$w^{t}: S^{N-1}\to T_{p_0}\M$ such that
\begin{equation}\label{eq:www}
S^{ {g}}_{t(1+{v^t})}(p_t)=
        \left\{\exp_{p_0}w^{{t}}(x)\,:\,x\in S^{N-1}\right\},\qquad \textrm{ and } \qquad  w^{t}(x)=\sum^{N}_{i=1}({t} x^i+ o({t}))E_i.
\end{equation}
Recall that
$$
S^{ {g}}_{t (1+v^t )}(p_t)=\left\{\exp_{p_t}\left( t(1+v^t)\sum^{N}_{i=1}x^i E_i^t \right)\,:\, x\in S^{N-1} \right\} .
$$
Let now consider the (well defined) map
$\Psi^{{t}}:=\exp^{-1}_{p_0}\circ\exp_{p_{{t}}}:T_{p_{{t}}}\M\longrightarrow
T_{p_{0}}\M$  and define
\begin{align*}
 F:[0,t_0)\times S^{N-1}\times T_{p_0}\M &\longrightarrow T_{p_0}\M\\
        ({t},x,w) &\longmapsto \Psi^{{t}}\biggl(t(1+
        v^{{t}}(x))\sum^{N}_{i=1}x^iE_{i}^t\biggl)-w.
 \end{align*}
 For a (fixed)  $x_0\in S^{N-1}$, we  have $F(0,x_0,0)=0$ and
 $D_{w}F(0,x_0,0)=-Id_{T_{p_0}\M}$. By the compactness of $S^{N-1}$, the implicit function theorem implies  that there
 exists ${t}_1 >0$  such that for all $ {t}\in (0,{t}_1)$ and  for all $x\in S^{N-1}$, there exits a
unique $w^{t}(x)\in T_{p_0}\M$ such that $F({t},x,w^{t}(x))=0$. That
is, for all $x\in S^{N-1}$
\begin{equation}\label{eq:www}
\exp_{p_0}(w^{t}(x))=\exp_{p_t}\left(t(1+
        {v^t}(x))\sum^{N}_{i=1}x^iE_{i}^t\right).
\end{equation}
In particular, we have $w(0,x)=0$ for all $x$ in $S^{N-1}$.
Differentiating  \eqref{eq:www} with respect to ${t}$, we get
\begin{equation}\label{d}
d(\exp_{p_0})_0\biggl(\frac{\partial w^
{t}}{\partial{t}}_{\mid_{{t}=0}}\biggl)=d(
\exp)_{p_0}\biggl(\frac{\partial
p_t}{\partial{t}}_{\mid_{{t}=0}}\biggl)+d(\exp_{p})_0\left(\sum^{N}_{i=1}x^iE_{i}\right).
\end{equation}
By assumption, we have  $\frac{\partial
p_t}{\partial{t}}_{\mid_{{t}=0}}=0$ and since
$$~~
 d(\exp_{p})_0=Id_{T_{p}\M},$$
we conclude that
$$\frac{\partial
w^{t}}{\partial{t}}_{\mid_{{t}=0}}=\sum^{N}_{i=1}x^iE_{i}
$$
~~and~~hence~
\be\label{eq:wteqtxo1}
w^{t}(x)=\sum^{N}_{i=1}  \left({t} x^i+O({t}^2 )\right)E_{i},
\ee
for all ${t}\in(0,{t}_1)$ and all $x\in S^{N-1}$ this proves the claim.\\

Observe that $|w^{{t}}(x)|_g\ne  0$ for ${t}>0$ small enough and thus we
can consider the map
$$
\alpha:(0,{t}_1)\times S^{N-1} \to S^{N-1}
$$
by
$$
 \a^i(t,x)=\frac{1}{|w^{t}(x)|_g } g(w^{t}(x), E_i ).
$$
It is clear from \eqref{eq:wteqtxo1} that
$$
\alpha({t},x)=\frac{x+O(t)}{|x+O(t)|}.
$$
The function $\alpha$ extends smoothly to ${t}=0$ with
$\alpha(0,\cdot)=Id_{S^{N-1}}$ and for  ${t}$ small enough
$\alpha({t},\cdot)$ is a diffeomorphism from $S^{N-1}$ into itself.
It is plain that for all $x\in S^{N-1}$
$$
w^t(x)=|w^t(x) |_g\,   \frac{w^t(x)}{|w^t(x) |_g}
$$
and thus for all     $y\in S^{N-1}$
\begin{equation}\label{pert2}
w^{t}(\alpha^{-1}({t},y))=|w^{t}(\alpha^{-1}({t},y))|_g\,
\sum_{i=1}^N y^i E_i.
\end{equation}
This together with \eqref{eq:www} imply that
\be \label{eq:BBw} S^{ {g}}_{t(1+{v^t})}(p_t) =\left\{\exp_{p_0}\biggl(|w^{t}(\alpha^{-1}({t},y))|_g\,\sum_{i=1}^N
y^i E_i\biggl)\,:\,~~ y\in S^{N-1}\right\}. \ee
 We have
\begin{equation}\label{diff}
|\alpha^{-1}({t},y)|^{2}=1~~\quad \textrm{ for~~all } ~~{t}\in(0,{t}_1)~~
\quad \textrm{ and }~~y\in S^{N-1}
\end{equation}
so that
\begin{equation}\label{difff}
\langle\partial_{{t}}\alpha^{-1}({t},y),\alpha^{-1}({t},y)\rangle=0.
\end{equation}
It then follows that
\begin{align*}
\partial_{{t}}(|w^{t}(\alpha^{-1}({t},y))|_g)
&=\frac{1}{|w^{t}(\alpha^{-1}({t},y))|_g}\langle
w^{t}(\alpha^{-1}({t},y)),\partial_{{t}}w^{t}+(d_xw^{t})(\partial_{{t}}\alpha^{-1})\rangle\\
&=\frac{1}{|\alpha^{-1}({t},y)+O(t)|}\langle
\alpha^{-1}({t},y)+O(t),\alpha^{-1}({t},y)+{t}\partial_{{t}}\alpha^{-1}+O(t)\rangle\\
&=\frac{1}{|\alpha^{-1}({t},y)+O(t)|}\langle
\alpha^{-1}({t},y)+O(t),\alpha^{-1}({t},y)+O(t)\rangle,
\end{align*}
where we have used \eqref{difff} to get the last line. Keeping in
mind that $\alpha (0,.)$ is the identity map, we obtain
$\partial_{{t}}(|w^{t}(\alpha^{-1}({t},y))|_g)_{\mid_{{t}=0}}=1$.
We conclude that  map
${t}\longmapsto|w^{t}(\alpha^{-1}({t},y))|_g$ is strictly
increasing with respect to ${t}\in (0,{t}_1)$ by decreasing
${t}_1>0$ if necessary. Therefore  thanks to \eqref{eq:BBw},
the family $\{S^{ {g}}_{{t}(1+{v^t})}(p_t),\,{t}\in
(0,{t}_0)\}$ constitutes a foliation and also setting
$ {\o}(t,y):= |w^{t}(\alpha^{-1}({t},y))|_g$, we finish the
proof of the proposition.\QED

\begin{Remark}\label{rem:foliation-PS}
An application of Proposition \ref{prop:folli} shows that the critical domains $\O_\e$ (in \eqref{eq:problemlam}) for the first eigenvalue of the Laplace-Beltrami operator  constructed by Pacard and Sicbaldi \cite{FPP}
constitutes also a local foliation of a neighborhood of the non-degenerate critical point $p_0$ of the scalar curvature. Indeed the improvement of the the distance between the center of their  extremal domains and $p_0$ was estimated by Sicbaldi and Dilay \cite{DS} which is of order $\e^2$.
\end{Remark}

\section{Proof of Theorem \ref{theo2}}\label{s:converse}
 Via the exponential map, we pull back the problem to   $ \R^{N}$. For this we consider the pull back metric
of $g$ under the map $\R^N\to \M,\;x \mapsto \exp_{p_0}\left(\e \sum_{i=1}^N x^i E_i\right)$, rescaled
with the factor $\frac{1}{\e^2}.$ Denoting this metric on $B_1$ by
$g_\e$, we then have, in Euclidean coordinates,
 \be \label{eq:dgvx00}
dv_{g_{\e}}(x):=\sqrt{|g_\e|}(x) = 1- O(\e^2).
 \ee
 Call $\Sig_\e=\de \O_\e \subset\R^N $  then it can be easily
verified that
$$
|\Sig_\e|_{g_\e}=|\Sig_\e|(1+O(\e^2))\quad |\O_\e|_{g_\e}=|\O_\e|(1+O(\e^2)).
$$
Integrate the first equality in \eqref{eq:problemz00} over $\O_\e$  to have
\begin{equation}\label{eq:PepoV}
|\Sig_\e|=\frac{N}{\e}|\O_\e|(1+o(1)).
\end{equation}
Now since $\O_\e\subset \d_\e B $ by \eqref{eq:ede}, we get that
$$
\frac{1}{\e}\leq\left( \frac{|B_1|}{|\O_\e| }\right)^\frac{1}{N}(1+o(1))
$$
and thus
$$
|\Sig_\e|\leq (1+o(1))N |\O_\e|^\frac{N-1}{N}|B_1|^\frac{1}{N}.
$$
We then conclude by the Euclidean  isoperimetric inequality that
\begin{equation}\label{eq:almost-min}
c_N |\O_\e|^\frac{N-1}{N}\leq |\Sig_\e|\leq
 (1+o(1))c_N |\O_\e|^\frac{N-1}{N},
\end{equation}
where $c_N=N|B|^\frac{1}{N}$ is the isoperimetric constant of  ${\R^{N}}$.
 In particular the sets $\O_\e$ are  almost minimizers for the isoperimetric problem.\\
Now consider the  real numbers $\rho_\e\to0$ defined as $|\O_\e|=|\rho_\e B_1|$.
Let $\Sig_\e'=\frac{1}{\rho_\e}\Sig_\e$ and $\O_\e'=\frac{1}{\rho_\e} \O_\e  $. Then    \eqref{eq:almost-min} yields
\begin{equation}\label{eq:almost-min-s}
 |S^{N-1}|\leq |\Sig_\e'|\leq
 (1+o(1))|S^{N-1}|,\quad |\O'_\e|=|B_1|.
\end{equation}
Using this and  \eqref{eq:PepoV} we get
\begin{equation}\label{eq:rovrho}
\frac{\rho_\e}{\e}=1+o(1)
\end{equation}
so that
\begin{equation}\label{eq:diam}
 \O_\e'\subset (1+o(1)) B_1.
\end{equation}
By compactness $\Sig_\e'$ converges weakly to $b+ S^{N-1}$  (see \cite{Maggi}) and also we have that   the symmetric distance $| \O_\e'\triangle (b+ B_1) |\to 0$  as $\e\to 0$, for some point $b\in \R^{N}$. Note that by \eqref{eq:diam}, $b=0$.
Letting $w_\e(x)=\rho_\e^{-2}u_\e(\rho_\e x)$, we have
\be\label{eq:op1P0}
\begin{cases}
\displaystyle -\D_{\ti{g}_\e} w_\e=1&\quad\textrm{ in } \O_\e'\vspace{3mm}\\
\displaystyle w_\e=0&\quad\textrm{ on } \de\O_\e'\vspace{3mm}\\
\displaystyle \ti{g}_\e(\n^{\ti{g}_\e} w_\e, \nu_\e')=-\frac{\e}{\rho_\e N}& \quad\textrm{ on } \de \O_\e',
\end{cases}
\ee
where $\ti{g}_\e(x)=g_\e(\rho_\e x)$.  It is also easy to see from \eqref{eq:uC2bound} that
\be\label{eq:v_rC2}
\|D^2 w_\e\|_{C(\ov{\O_\e'})}\leq c  .
\ee
We let $d_\e(x)=\textrm{dist}(\de\O_\e',x)$ be the distance
function of $\Sig_\e'$. Given $x\in\O_\e'$  near $\de \O_\e'$ then it can
be written uniquely as $x=\s_x-d_\e(x)\,\nu_\e'(\s_x)$, where
$\s_x$ is the projection of $x$ on $\Sig_\e'$. This defines coordinates $(t,\s)\mapsto x=\s-t\nu_\e'(\s)$.
Recall the decomposition of the Laplace-Beltrami operator in the coordinates $(t,\s)$:
$$
\D_{\ti{g}_\e}=\frac{\de^2}{\de t^2}+H^t_\e \frac{\de}{\de t}+\D_{\Sig^t_\e},
$$
where $H_\e^t$ is the mean curvature of  the hypersurface $\Sig^t_\e=\{x\in \O'_\e\,:\, d_\e=t \}$ with respect to the metric $ \ti{g}_\e$ and $\D_{\Sig^t_\e}$ is the Laplace-Beltrami on   $\Sig^t_\e$. We also observe that
$$
\frac{\de w_\e}{\de t}=|\n w_\e|_{\ti{g}_\e}=-\ti{g}_\e(\n_{\ti{g}_\e} w_\e, \nu_\e') \quad \textrm{ on }\de\O_\e'.
$$
Thanks to \eqref{eq:v_rC2} and the second equation in  \eqref{eq:problemz00}, we conclude that
$$
H_\e^0= \frac{1-\frac{\de^2w_\e}{\de t^2} }{\frac{\de w_\e}{\de t}} \qquad \textrm{ on } \de \O_\e'.
$$
Therefore
$$
|H_\e^0|\leq Const.\qquad \textrm{ on } \de \O_\e'.
$$
Since $\ti{g}_\e$ is nearly Euclidean, the mean curvature of $\de\O_\e'$, with respect to the Euclidean metric, is uniformly bounded with respect to $\e$. Hence
by \cite{Nard} (see also \cite{FM-DJ}) the hypersurface $\Sig_\e'$ converges smoothly to $S^{N-1}$ and  there exists a
 function  $\ti{v}^\e\in C^{2,\a}(S^{N-1})$ with
 $\|\ti{v}^\e\|_{C^{2,\a}(S^{N-1})}\to 0$ as $\e\to0$ and  such that
$$
\Sig'_\e=(1+\ti{v}^\e)  S^{N-1}.
$$
We therefore conclude from \eqref{eq:rovrho} that
$$
\Sig_\e=\rho_\e(1+\ti{v}^\e)  S^{N-1}=\e(1+ {v}^\e) S^{N-1}
$$
and of course $\| {v}^\e\|_{C^{2,\a}(S^{N-1})}\to 0$ as $\e\to0$.
Hence we   get  $\O_\e=B^g_{\e(1+  {v}^\e)}(p_0)$ so that the uniqueness of Proposition \ref{Propo3.2} and a scaling argument yield
$$
u_\e=\e^2\bar\phi(p,\e,v^\e_0,\bar{v}^\e).
$$
Since, by assumption,
\be \label{eq:gnuenue}
\e^{-1}{g}(\nabla_{ {g}}u_\e, {\nu_\e})|_{\de B^{  g}_{\e(1+ {v}^\e)}}=   \hat{g}(\nabla_{\hat{g}} \hat{u}_\e ,\hat{\nu}_\e)|_{\de B_1}  =-\frac{1}{N},
\ee
the uniqueness of
  Proposition \ref{Pro3.5}  implies that
\be\label{eq:Pi1ve0}
\Pi_1{v}^\e=0
\ee
 provided $\e$ is small.
We now compute the normal derivative of $u_\e$ by using similar arguments as in the proof of Lemma \ref{Lem3.3}. It follows that
\begin{align*}
\e^{-1}  {g}(\nabla_{ {g}}u_\e, {\nu_\e})|_{\de B^{  g}_{\e(1+ {v}^\e)}} =-\frac{1}{N} +\frac{1}{N}\mbL(  {v}^\e )  +  (\partial_{\nu}\psi_\e)_{|_{\partial
B_1}}+(\partial_{\nu}\G_{\e,v^\e})_{|_{\partial
B_1}}   +P_{\e}^1(v^\e),
\end{align*}
where $\G_{\e,v}$ satisfies \eqref{eq:Gev} and  the
function   $\psi_{\e}$  satisfies
\begin{align}\label{eq:eqo0}
\begin{cases}
\displaystyle -\Delta\psi_\e=\frac{\e^2}{3N}Ric_{p_0}(X,X) -\frac{\e^3}{4N} g(\n_XR_{p_0}(E_i,X)E_i,X)\\
\qquad\qquad \qquad\displaystyle +\frac{\e^3}{6N}g(\n_XR_{p_0}(E_i,X)X,E_i)+O_{p_0}(\e^4)& \quad \textrm{ in }  B_1 \vspace{3mm}\\
\displaystyle\psi_\e=0&\quad \textrm{ on } \de B_1.
\end{cases}
\end{align}
Thanks to \eqref{eq:gnuenue}, we  have
\be\label{eq:mbLve}
\frac{1}{N} \mbL(  {v}^\e )=
 - (\partial_{\nu}\psi_\e)_{|_{\partial
B_1}}+(\partial_{\nu}\G_{\e,v^\e})_{|_{\partial
B_1}}   +P_{\e}^1(v^\e).
\ee
From \eqref{eq:Gev}, we see immediately from elliptic regularity theory that
$$
\|\G_{\e,v^\e}\|_{C^{2,\a}(S^{N-1})}\leq C\e^4+ C\e^2\|v^\e\|_{C^{2,\a}(S^{N-1})}+C\|v^\e\|_{C^{2,\a}(S^{N-1})}^2.
$$
Recalling   \eqref{eq:Pi1ve0}, we then   apply  Proposition \ref{Propo3.4} in \eqref{eq:mbLve} to have
$$
\|v^\e\|_{C^{2,\a}(S^{N-1})}\leq C\e^2+C\e^2 \|v^\e\|_{C^{2,\a}(S^{N-1})}+C\|v^\e\|_{C^{2,\a}(S^{N-1})}^2.
$$
This implies that
$$
\|v^\e\|_{C^{2,\a}(S^{N-1})}\leq C\e^2.
$$
We then conclude that
$$
 \mbL(  {v}^\e )=
 -  N   (\partial_{\nu}\psi_\e)_{|_{\partial
B_1}}  +  O(\e^4)
$$
Now we multiply this equation  by $x^i$, integrate by parts over $B_1$,  use \eqref{eq:eqo0} together with Bianchi's  identity to get
$$
\n^i_g S_g(p_0)=0.
$$
\QED

\section{Appendix: }
As mentioned in the first section, the torsional rigidity of the rod
$\O\times \R$  is proportional to the inverse of
\begin{equation}\label{eq:propp99}
J(\Omega):=\inf\biggl\{\int_{\Omega}|\n
u|^{2}_g\,\textrm{dvol}_g:~~\int_{\Omega}u\,\textrm{dvol}_g=1,~u\in
H^{1}_{0}(\Omega)\biggl\}.
\end{equation}
 In particular minimizing $\O\mapsto J(\O)$ is equivalent to maximizing the torsion rigidity and therefore Serrin's
 result states that balls maximize the torsion rigidity as it can be also derived from the Faber-Krahn inequality.\\
In this appendix we consider the  isochoric profile for the torsion problem defined as
 \begin{equation}\label{eqproo}
 \calT_\M(v,g):=\inf_{ \Omega\subset
 \M,|\Omega|_g=v}J(\Omega),
 \end{equation}
where here and in the following, we assume without further mention that only
regular bounded domains $\O\subset\M$ are considered. In particular thanks to the Faber-Krahn inequality
$$
\calT_{\R^N}(v)=J(B_1)\biggl(\frac{|B_1|}{v}\biggl)^{-\frac{N+2}{N}}.
$$
Similarly in the space of constant sectional curvatures, balls minimize $J$, see \cite{Mc}.
Isochoric comparison for   $\calT$ has been studied recently in the papers \cite{GLL}, \cite{Xiao}.
 Here we deal with local asymptotics of this profile as $v\to0$.  This also leads to
isochoric comparison in terms of scalar curvature.

 In the recent years, several works have been devoted to the  Taylor expansion of isoperimetric
 and ischoric profile for some geometric quantities such as  the (relative) perimeter functional,
 Cheeger constants,   Dirichlet eigenvalue, second Neuman eigenvalue, etc. We refer the papers
  \cite{Druet-FK,OD,Nard, mmf,mmf1,  BM, Bayle-Rosales, FW, {Chavel-1}}.
   We should mention that the argument in this section will follow closely
    Druet \cite{OD} where he studied the expansion of the Faber-Krahn profile.
    The main result of this section is contained in the following
 \begin{Theorem}\label{th:Isoco-profile}
Let $(\M,g)$ be a compact Riemannian manifold of dimension $N\geq 2$
. As $v\to 0$, we have
\begin{align*}
 \calT_\M(v,g)&=
  \biggl[1-\frac{N+6}{6N(N+4)}\biggl(\frac{v}{|B_1|}\biggl)^{\frac{2}{N}}\max_{\M}S_g+O(v^{\frac{3}{N}})\biggl]\calT_{\R^N}(v),
\end{align*}
 where  $S_g$ is the scalar
curvature of $(\M,g)$.
\end{Theorem}
\proof The first step of the proof is to derive the expansion of
$J(B^g_\e(p))$ as $\e\to 0$. Once this is done we then obtain an
upper bound for $\calT_\M(v,g)  $ as $v\to 0$. The second step
consists in using  the asymptotic profile of the isoperimetric
profile for the perimeter functional obtained by Druet in \cite{OD}
together with the Faber Krahn inequality on the space of constant
sectional curvatures.
This later step follows exactly Druet \cite{Druet-FK}. Therefore we will only give the proof of the first step.\\

\noindent
\textbf{Claim:} As $v\to 0$, we have
\be\label{eq:claTvgle}
 \calT_\M(v,g)\leq
  \biggl[1-\frac{N+6}{6N(N+4)}\biggl(\frac{v}{|B_1|}\biggl)^{\frac{2}{N}}\max_{\M}S_g+O(v^{\frac{3}{N}})\biggl]\calT_{\R^N}(v).
\ee To see this we determine the Taylor expansion of $J(B^g_\e(p))$
as $\e\to 0$. Recall that   $J(B^g_\e(p))$ is the  Dirichlet energy
 in the ball $B^g_\e(p)$ and $u_{\e}$ the corresponding minimizer,  that is
 \begin{equation}\label{eq:pl}
 \begin{cases}
 -\Delta_{g}u_{\e}=J(B^g_\e(p))& \quad\textrm{ in } ~B^g_\e(p)\\
 ~u_{\e}=0&~\quad\textrm{ on } ~\partial
 B^g_\e(p)\\
 \displaystyle \int_{ B^g_\e(p)}u_{\e}\,\textrm{dvol}_g=1.
 \end{cases}
 \end{equation}
More precisely, we have that
\begin{equation}\label{eq:pm}
J(B^g_\e(p))\leq\int_{B^g_\e(p)}|\nabla_{g}u|^{2}\textrm{dvol}_{g}
\end{equation}
 for all $u\in H^{1}_{0}(B^g_\e(p))$ such that
 $\int_{B^g_\e(p)}u_{\e}\textrm{dvol}_{g}=1$.
  Via the exponential map, we pull back the problem to the unit
ball $B_1 \subset \R^{N}$. For this we consider the pull back metric
of $g$ under the map $B_1\to \M,\;x \mapsto \exp_p(\e x)$, rescaled
with the factor $\frac{1}{\e^2}.$ Denoting this metric on $B_1$ by
$g_\e$, we then have, in Euclidean coordinates, \be \label{eq:dgvx}
\textrm{dvol}_{g_{\e}}(x)=\sqrt{|g_\e|}(x) = 1-\frac{\e^2}{6}
Ric_{p}(X,X) +O(\e^3) \ee
 for $x \in \overline B_1$ by Proposition \ref{prop:expgij}.
 We consider the function $\varphi_{\e}(x):=\e^{N}u_{\e}(\e x)$ and we recall in \eqref{exp} the expansion of the scaled metric
 $\ti{g}_{\e}(x)=g(\e x)$ for $x\in B_1$. From
\eqref{eq:pl}, we get
 \begin{equation}\label{eq:pn}
 \begin{cases}
-\Delta_{g_{\e}}\varphi_{\e}=J(B^g_\e(p))\e^{N+2}& \quad \textrm{ in } ~B_1 \vspace{3mm}\\
\varphi_{\e} =0 &\quad \textrm{ on } ~\de B_1\vspace{3mm} \\
\displaystyle\int_{B_1}\varphi_{\e}(x)\textrm{dvol}_{g_{\e}}=1.
\end{cases}
\end{equation}
The functions $\varphi_{\e}$ are positive in $B_1$ and equal
to $0$ on the boundary. Thank to \eqref{eq:pm}, we obtain
\begin{equation}\label{eq:po}
J(B^g_\e(p))\e^{N+2}\leq\int_{B_1}|\nabla_{g_{\e}}
\varphi_\e|^{2}\textrm{dvol}_{g_{\e}}
\end{equation}
for all $u\in H^{1}_{0}(B_1)$ such that
 $\int_{B_1}u\textrm{dvol}{g_{\e}}=1$.
 Since the metric $g_{\e}\longrightarrow g_{0}$ as
 $\e\longrightarrow0$ this immediately implies $\limsup_{\e\rightarrow0}J_{\e}\e^{N+2}\leq J_1$,
 where $J_1=J(B_1)$. Using equation \eqref{eq:pn} and  regularity results, the sequence $(\varphi_{\e})$ is
 uniformly bounded in $C^{2}(B_1)$, and we can write
$\int_{B_1}\varphi_{\e}=1+O_{p}(\e)$ and $\int_{B_1}|\nabla
\varphi_{\e}|^{2}=\e^{N+2}J(B^g_\e(p))+O_{p}(\e)$. This implies that
$\liminf_{\e\rightarrow0}\e^{N+2}J(B^g_\e(p))\geq J_1$. So we have
proved that $J(B^g_\e(p))\e^{N+2}\longrightarrow J_1$ as
$\e\longrightarrow0$. Since $(\varphi_{\e})$ is
 uniformly bounded in $C^{2}(B_1)$ and any subsequence has to converge
 to the (unique) solution of the limit equation
 $-\Delta \varphi=J_{1}$ in $B_1$ with $\int_{B_1}\varphi=1$ and $\varphi\geq0$, we deduce that
 $\varphi_{\e}\longrightarrow \varphi$ in $C^{1} (\ov{B_1}) $ as $\e\longrightarrow0$. We multiply
 \eqref{eq:pn} by $\varphi$ and we get after integrating by parts,
 $$J(B^g_\e(p))\e^{N+2}=-\int_{B_1}\varphi_{\e}\Delta_{g_{\e}}\varphi \textrm{dvol}_{g_{\e}} .$$ By Lemma \ref{Lem2}  we have
\begin{align*}
\Delta_{g_{\e}}\varphi&=\Delta\varphi-\frac{1}{3}\sum^{N}_{k,l,i,j=1}R_{ikjl}x^{k}x^{l}\e^{2}\partial^{2}_{ij}\varphi+\frac{2}{3}\sum^{N}_{i,j,k=1}R_{ijik}x^{k}\e^{2}\partial_j\varphi+O_p(\e^{3})\\
&=-J_1-\frac{1}{3}\sum^{N}_{k,l,i,j=1}R_{ikjl}x^{k}x^{l}\e^{2}\partial^{2}_{ij}\varphi+\frac{2}{3}\sum^{N}_{i,j,k=1}R_{ijik}x^{k}\e^{2}\partial_j\varphi+O_p(\e^{3})
\end{align*}
and thus
\begin{align*}
J(B^g_\e(p))\e^{N+2}-J_1\int_{B_1}\varphi_{\e}\textrm{dvol}_{g_{\e}}  &=\frac{1}{3}\e^{2}\sum^{N}_{k,l,i,j=1}R_{ikjl}\int_{B_1}x^{k}x^{l}\varphi_{\e}\partial^{2}_{ij}\varphi \\
 &-\frac{2}{3}\e^{2}\sum^{N}_{i,j,k=1}R_{ijik}\int_{B_1}x^{k}\varphi_{\e}\partial_j\varphi  +O_p(\e^{3}).
\end{align*}
Thanks to \eqref{eq:pn} and the convergence of $\varphi_{\e}$
to $\varphi$, a straightforward computation using  also
\eqref{eq psi1} yields
\begin{align*}
\sum^{N}_{k,l,i,j=1}R_{ikjl}\int_{B_1}x^{k}x^{l}\varphi\partial^{2}_{ij}\varphi
&=-\frac{J^{2}_1}{N}\sum^{N}_{k,l,i=1}R_{ikil}\int_{B_1}x^{k}x^{l}\phi_0
\\
&=\frac{J_1}{N(N+4)}S_g(p)=\sum^{N}_{i,j,k=1}R_{ijik}\int_{B_1}x^{k}\varphi\partial_j\varphi
\end{align*}
 and thus  we get
 \be \label{eq:expanJBe}
 J(B^g_\e(p))\e^{N+2}
=J_1\biggl(1-\frac{S_g(p)}{3N(N+4)}\e^{2}+O_p(\e^{3})\biggl),
\ee where $\varphi$ is given by
$\varphi=\phi_{0}/||\phi_{0}||_{L^{1}(B_1)}$ and $\phi_0$ is the
unique solution of \eqref{eq:ropp2}. Next recall the expansion of
volume of geodesic balls which can be deduce from \eqref{eq:dgvx}:
$$
\left|B^g_\e(p) \right|_g= \e^{N}\,{\left|B_1
\right|}\,\left(1-\frac{1}{6(N+2)}\,\e^2{S}_g(p)+O(\e^3)  \right).
$$
This together with  \eqref{eq:expanJBe} implies that
$$
J(B^g_\e(p))=
\biggl[1-\frac{N+6}{6N(N+4)}\biggl(\frac{v}{|B_1|}\biggl)^{\frac{2}{N}}
S_g(p)+O(v^{\frac{3}{N}})\biggl]\calT_{\R^N}(v),
$$
where $v= \left|B^g_\e(p) \right|_g$. This then proves the claim as $p$ is arbitrary.\\

As said above, the reverse inequality of \eqref{eq:claTvgle}
follows step-by-step \cite{OD} so we skip the details.

\QED


\begin{thebibliography}
\footnotesize

\bibitem{Alexandrov} A. D. Alexandrov,  {Uniqueness Theorem
for surfaces in large I}, Vestnik Leningrad Univ. Math. 11 (1956),
5-17.
%
\bibitem{AB} A. Ambrosetti and M. Badiale, Variational perturbative methods and bifurcation of bound states from the essential spectrum, Proc. Roy. Soc. Edinburgh Sect. A 128, (1998), 1131-1161.
 %
\bibitem{AM} A. Ambrosetti and A. Malchiodi, Perturbation Methods and Semilinear Elliptic Problems on $R^n$. Progress in Mathematics, Birkh\"auser Verlag, Basel-Boston-Berlin (2005).
%
%
\bibitem{BM} P. Berard  and D. Meyer, In\'egalit\'es isop\'erim\'etriques et applications.
Ann. Sci. \'Ec. Norm. Sup\'er., IV. S\'er. 15, (1982) 513-541.

\bibitem{Bayle-Rosales} V. Bayle  and  C.  Rosales, Some isoperimetric comparison theorems
for convex bodies in Riemannian manifolds, Indiana Univ. Math. J. 54 (2005),
1371-1394.
%
%
\bibitem{Chavel-1} I.  Chavel, Eigenvalues in Riemannian geometry, Academic Press, 1984.
%
\bibitem{CH} M. Choulli and A. Henrot,  Use of the Domain Derivative to Prove Symmetry Results in Partial Differential Equations.
Mathematische Nachrichten, Volume 192, Issue 1, pages 91-103, 1998.

\bibitem{DS} E. Daley and P. Sicbaldi {Extremal domains for the first
eigenvalue of the Laplace Beltrami operator in a general compact
riemannian manifold} arXiv:1302.4221v1.

\bibitem{dPPW} M. del Pino, F. Pacard and J. Wei, Serrin's Overdetermined problem and constant mean curvature surfaces.
Preprint.  http://arxiv.org/abs/1310.4528v1


\bibitem{doCarmo} M. do Carmo, Riemannian Geometry. Boston: Birkhauser.  (1992).

\bibitem{Druet-FK} O. Druet,  Asymptotic expansion of the Faber-Krahn profile of
a compact Riemannian manifold. C. R. Math. Acad. Sci. Paris  346
(2008),  no. 21-22, 1163-1167.
%

\bibitem{OD} O. Druet, Sharp local isoperimetric inequalities involving the scalar curvature, Proceedings of
 the American Mathematical Society, 130, 8, 2351-2361, (2002).

\bibitem{ELIL}A. El Soufi and S. Ilias, Domain deformations and eigenvalues of
the Dirichlet Laplacian in Riemannian manifold, Illinois Journal of
Mathematics 51 (2007) 645-666.




\bibitem{mmf}M. M. Fall, Area-minimizing  sets in Riemannian manifolds with boundary constrained to small volume.  Pacific J. Math. 244 (2010), no. 2, 235--260.
%
\bibitem{mmf1}M. M. Fall, Some local eigenvalue estimates  involving   curvatures.  Calc. Var. Partial Differential Equations  36  (2009),  no. 3, 437-451.

\bibitem{mmfcm1}  M. M. Fall and C. Mercuri, Minimal disc-type surfaces embedded in a perturbed cylinder. Differential Integral Equations 22 (2009), no. 11-12, 1115--1124.

\bibitem{mmfcm} M. M. Fall  and C. Mercui, Foliations of small tubes in Riemannian manifolds by capillary minimal discs,  Nonlinear Analysis TMA.  (12) Vol. 70 (2009) 4422-4440.

\bibitem{mmm} M. M. Fall  and F.   Mahmoudi, Hyper-surfaces with free boundary and large constant mean curvature: concentration
 along sub-manifolds. Ann. Scuola Norm. Sup. Pisa Cl. Sci. (5) Vol. VII (2008), 1-40.

\bibitem{FW} M. M. Fall and T.  Weth,   Sharp local upper bound of the first non-zero Neumann eigenvalue in Riemannian manifolds. To appear in Calc. Var. Partial Differential Equations. http://arxiv.org/abs/1110.4770.

\bibitem{FK}  A. Farina and B. Kawohl,   {Remarks on overdetermined
boundary value problem } Calc. Var. Partial Differential Equations
31 (2008), n0 3. 351-357.
%


%

\bibitem{GLL} S. Gallot, A. Loi and C. Lucio, Maximizing torsional rigidity on Riemannian manifolds. http://arxiv.org/abs/1309.7796.

\bibitem{He} D. Henry, Perturbation of the Boundary in Boundary-Value Problems
of Partial Differential Equations. Cambridge University Press.  2005.

\bibitem{FM-DJ} D. L. Johnson D.L. and F. Morgan,  Some sharp isoperimetric
theorems for riemannian manifolds. Indiana Univ. Math. J., (2000) 49 (2).

 \bibitem{Kap} N.  Kapouleas, Compact constant mean curvature surfaces in Euclidean
three-space. J. Differ. Geom. 33, No.3, (1991), 683-715.

 \bibitem{Landau} L. D. Landau and E.M. Lifshitz, Theory of Elasticity, Course of Theoretical Physics
7, Pergamon Press, London, 1959.

%
\bibitem{Laurin} P. Laurin,  Concentration of CMC surfaces in a Riemannian manifold,  Int. Math. Res. Not. IMRN 2012, no. 24, 5585-5649.
%

\bibitem{Lind} F. Lindemann,  \"{U}ber die Zahl $\pi$. Math. Ann. 20, 213-225, 1882.

\bibitem{mmp} F. Mahmoudi, R.  Mazzeo and F. Pacard, Constant mean
curvature hypersurfaces condensing along a submanifold, Geom. funct.
anal. Vol. 16 (2006) 924-958.
%
\bibitem{mp} R. Mazzeo, F. Pacard, Foliations by constant mean
curvature tubes,  Comm. Anal. Geom.  13  (2005),  no. 4, 633-670.
%

\bibitem{Maggi} F. Maggi, Sets of finite perimeter and geometric variational problems: an introduction to Geometric
Measure Theory, Cambridge Studies in Advanced Mathematics no. 135, Cambridge University Press,
2012.
%

 \bibitem{malm}A.  Malchiodi  and M.  Montenegro, Boundary concentration
phenomena for a singularly perturbed elliptic problem, Comm. Pure
Appl. Math. 55 (2002), no. 12, 1507-1568.
%
\bibitem{Mc} P. Mcdonald,  Isoperimetric Conditions, Poisson Problems, and
Diffusions in Riemannian Manifolds. Potential Analysis 16: 115-138, 2002.
%
\bibitem{MiPi} A. M. Micheletti, A. Pistoia, Generic properties of critical points of the scalar curvature for
a Riemannian manifold, Proc. Amer. Math. Soc. 138 (2010), no. 9, 3277-3284.

\bibitem{MS}F. Morabito and P. Sicbaldi, Delauney type domains for an overdetermined elliptic problem in $S^n\times \R$ and $\mathbb{H}\times \R$. Preprint. http://arxiv.org/abs/1305.6516

 \bibitem{Nard} S. Nardulli,  R\'egularit\'e des solutions du probl\'eme
isop\'erim\'etrique proches de vari\'et\'es lisses. Preprint Universit\'e
de Paris Sud, Octobre 2006.
%
%
%

\bibitem{FPP} F. Pacard and P. Sicbaldi, Extremal domains for the first eigenv
alue of the Laplace-Beltrami operator. Ann. Inst. Fourier (Grenoble)
59 (2009), no. 2, 515-542.

\bibitem{PX} F. Pacard, X. Xu, Constant mean curvature sphere in Riemannian
manifolds. Manuscripta Math. 128 (2009), no. 3, 275-295.
%
%
\bibitem{RaSa} S. Raulot and  A. Savo,  On the spectrum of the Dirichlet-to-Neumann operator acting on forms of a Euclidean domain. J. Geom. Phys. 77 (2014) 1-12.
%
\bibitem{ScSi} F. Schlenk and  P. Sicbaldi, Bifurcating extremal domains for the
first eigenvalue of the Laplacian. Adv. Math. 229 (2012) 602-632.
%
%



\bibitem{STY} R. Schoen and S.T. Yau, Lectures on Differential Geometry,
International Press (1994).

%
\bibitem{Serrin} J. Serrin, A Symmetry Theorem in Potential Theory. Arch. Rational
Mech. Anal. 43 (1971), 304-318.

\bibitem{Sic} P. Sicbaldi, New extremal domains for the first eigenvalue
of the Laplacian in flat tori. Calc. Var. (2010) 37:329-344.
%
\bibitem{Sun} T. Sun, A note on constant geodesic curvature curves on surfaces. Ann. Inst. H.
Poincaré Anal. Non Lin\'eaire, 26(5) :1569-1584, 2009.
%
\bibitem{Wein} H. F. Weinberger,
An isoperimetric inequality for the $N$-dimensional free membrane problem. J. Rational Mech. Anal. 5 (1956), 633-636.
%

\bibitem{TJW} T.J. Willmore, Riemannian Geometry, Oxford Univ. Press. NY.
(1993).



\bibitem{Xiao} J. Xiao,  Isoperimetry for semilinear torsion problems in
Riemannian two-manifolds. Advances in Mathematics 229 (2012) 2379-2404.

%
\bibitem{Ye}  R. Ye  {Foliation by constant mean curvature spheres
}. Pacific J. Math. 147 (1991), no. 2, 381-396.
\end{thebibliography}
\end{document}